\documentclass[11pt,a4paper]{article}
\usepackage[utf8]{inputenc}
\usepackage{amsmath}
\usepackage{amsfonts}
\usepackage{color}
\usepackage{amssymb}
\usepackage{mathrsfs}
\usepackage{esint}
\usepackage[colorinlistoftodos]{todonotes}
\usepackage[colorlinks=true]{hyperref}

\hypersetup{
	colorlinks=blue,%
	citecolor=red,%
	filecolor=black,%
	linkcolor=blue,%
	urlcolor=blue
}
\usepackage[margin=1.3in]{geometry}

\usepackage[amsmath,thmmarks,hyperref]{ntheorem}
{
	\theoremstyle{nonumberplain}
	\theoremheaderfont{\bfseries}
	\theorembodyfont{\normalfont}
	\theoremsymbol{\mbox{$\Box$}}
	\newtheorem{pf}{Proof.}
}

\numberwithin{equation}{section}
\def\R{\mathbb{R}}
\def\B{\mathbb{B}}
\def\S{\mathbb{S}}

\def\e{\epsilon}
\newcommand{\ud}{\mathrm{d}}

\newtheorem{thm}{Theorem}[section]
\newtheorem{definition}{Definition}[section]

\newtheorem{lem}{Lemma}[section]
\newtheorem{rem}{Remark}[section]

\newtheorem{pro}{\indent Proposition}[section]

\newtheorem*{convention}{\textbf{Convention.}}
\newtheorem*{thm A}{Theorem A}
\newtheorem{cor}{Corollary}[section]
\usepackage{upgreek}
\usepackage{graphicx}

\makeatletter

\newdimen\bibspace
\renewenvironment{thebibliography}[1]{%
	\section*{\refname 
		\@mkboth{\MakeUppercase\refname}{\MakeUppercase\refname}}%
	\list{\@biblabel{\@arabic\c@enumiv}}%
	{\settowidth\labelwidth{\@biblabel{#1}}%
		\leftmargin\labelwidth
		\advance\leftmargin\labelsep
		\itemsep\bibspace
		\parsep\z@skip     %
		\@openbib@code
		\usecounter{enumiv}%
		\let\p@enumiv\@empty
		\renewcommand\theenumiv{\@arabic\c@enumiv}}%
	\sloppy\clubpenalty4000\widowpenalty4000%
	\sfcode`\.\@m}
{\def\@noitemerr
	{\@latex@warning{Empty `thebibliography' environment}}%
	\endlist}

\makeatother

\makeatletter

\allowdisplaybreaks

\allowdisplaybreaks

\allowdisplaybreaks

\begin{document}
	\title{\bf  A simple proof of reverse Sobolev inequalities on the sphere and  Sobolev trace inequalities on the unit ball}
	\date{}

\author{Runmin Gong\thanks{R. Gong, School of Mathematics and Statistics, Wuhan University, Wuhan, 430072, People's Republic of China, Email:2022302011151@whu.edu.cn. }, Qiaohua  Yang\thanks{Q. Yang, School of Mathematics and Statistics, Wuhan University, Wuhan, 430072, People's Republic of China, Email: qhyang.math@whu.edu.cn.}, Shihong Zhang\thanks{ S. Zhang, School of Mathematics \& IMS, Nanjing University, Nanjing 210093, People's Republic of China, Email:dg21210019@smail.nju.edu.cn.}}

	



	\renewcommand{\thefootnote}{\fnsymbol{footnote}}
	\maketitle

	{\noindent\small{\bf Abstract: Frank et al. (J. Funct. Anal., 2022) stated that there is no relation  between the reversed Hardy-Littlewood-Sobolev (HLS) inequalities and reverse Sobolev inequalities. However, we demonstrate that reverse Sobolev inequalities of order
			 $\gamma\in(\frac{n}{2},\frac{n}{2}+1)$ on the $n$-sphere can be readily derived from the reversed HLS inequalities. For the case
			  $\gamma\in(\frac{n}{2}+1,\frac{n}{2}+2)$,  we present a simple proof of reverse Sobolev inequalities by using the center of mass condition
introduced by Hang.  In addition, applying this approach, we  establish the quantitative stability of reverse Sobolev inequalities of order  $\gamma\in(\frac{n}{2}+1,\frac{n}{2}+2)$  with explicit lower bounds.  Finally, by using conformally covariant boundary operators and reverse Sobolev inequalities, we derive  Sobolev trace inequalities on the unit ball.}

		\medskip
		
		{{\bf $\mathbf{2020}$ MSC:} 43A90, 46E35 (58J32, 31B30)}
		
		
		
		\section{Introduction}	
		
	Conformally invariant operators have long been a cornerstone in the study of conformal geometry. The most general such operators on  $\mathbb{S}^{n}=\{x\in\mathbb{R}^{n+1}:|x|=1\}$ are the GJMS operators $P_{2\gamma}$, whose general expression is given by the following:	
		\begin{align*}
		P_{2\gamma}=\frac{\Gamma\left(B+\frac{1}{2}+\gamma\right)}{\Gamma\left(B+\frac{1}{2}-\gamma\right)}, \qquad \mathrm{where}\qquad B=\sqrt{-\Delta_{\S^n}+\frac{(n-1)^2}{4}}.
		\end{align*}
	This implies that for any $ f \in C^{\infty}(\mathbb{S}^n)$ with the spherical expansion
	$$ f = \sum_{l=0}^{\infty} Y_l, \,\,\quad Y_l \in \mathscr{H}_l, $$
	where $\mathscr{H}_l $ denotes the space of spherical harmonics of degree  $ l$, one has
		\begin{align}\label{1.1}
			P_{2\gamma}f=\sum_{l=0}^{+\infty}\frac{\Gamma\left(l+\frac{n}{2}+\gamma\right)}{\Gamma\left(l+\frac{n}{2}-\gamma\right)}Y_{l}.
		\end{align}

	The Sobolev inequalities involving these GJMS operators have been developed over many years. The simplest case occurs when $\gamma = 1 $, which asserts that for any $ f \in H^1(\mathbb{S}^n)$, the following inequality holds (see \cite{Aubin, Talenti})
		\begin{align}\label{Intr 2rd Sobolev}
			\int_{\S^n}\left(|\nabla_{\S^n}f|^2+\frac{n(n-2)}{4}f^2\right)   \ud V_{{\S^n}}\geq \frac{n(n-2)}{4}|\S^n|^{\frac{2}{n}}\|f\|^2_{L^{\frac{2n}{n-2}}(\S^n)},
		\end{align}
where $|\mathbb{S}^n|$ denotes the volume of $\mathbb{S}^n$. This inequality is equivalent to the classical Sobolev inequality in $ \mathbb{R}^n $, namely, for any $ u \in H^1(\mathbb{R}^n) $,
		\begin{align*}
			\int_{\R^n}|\nabla u(x)|^2dx\geq \frac{n(n-2)}{4}|\S^n|^{\frac{2}{n}}\|u\|^2_{L^{\frac{2n}{n-2}}(\R^n)}.
		\end{align*}
	For general $ \gamma \in \left(0, \frac{n}{2}\right) $, the sharp Sobolev inequalities are established in dual form by Lieb in \cite{Lieb} (see also Frank and Lieb \cite{Frank&Lieb}). The case $ \gamma = \frac{n}{2} $  is known as the Onofri-Beckner inequality (see Onofri \cite{on} for $n=2$ and Beckner \cite{Beckner} for $n\geq 3$), which can be viewed as the limiting form of the Sobolev inequalities as $ \gamma \to \frac{n}{2} $ (see \cite{Chang&Wang, Xiong}). We  present these results as follows:
		\begin{thm A*}
If $0<\gamma<\frac{n}{2}$,  then we have, for
			any  $f\in H^{\gamma}(\S^n)$,
			\begin{align}\label{Intro General ineu}
			\int_{\S^n}fP_{2\gamma}f   \ud V_{{\S^n}}\geq \frac{\Gamma\left(\frac{n}{2}+\gamma\right)}{\Gamma\left(\frac{n}{2}-\gamma\right)}|\S^n|^{\frac{2\gamma}{n}}\|f\|^2_{L^{\frac{2n}{n-2\gamma}}(\S^n)}.
			\end{align}
			The equality holds if and only if there exists a conformal transformation $\varphi:\mathbb{S}^{n}\rightarrow\mathbb{S}^{n}$ such that
$$f=c(\det\ud \varphi)^{\frac{n-2\gamma}{2n}},\;c\in \mathbb{R}.$$
 If $\gamma=\frac{n}{2}$, then we have, for any $f\in H^{\frac{n}{2}}(\S^n)$,
			\begin{align}\label{Intro Beckner}
			\frac{1}{2n!}\fint_{\S^n} fP_{n}f  \ud V_{{\S^n}}\geq \log\left(\fint_{\S^n} e^{f-\bar{f}}  \ud V_{{\S^n}}\right),
			\end{align}
where $$\fint_{\S^n} f  \ud V_{{\S^n}}=\frac{1}{|\S^{n}|}\int_{\S^n} f  \ud V_{{\S^n}}\;\;\;
 \textrm{and} \;\;\;\bar{f}=\fint_{\S^n} f  \ud V_{{\S^n}}.$$
The equality holds if and only if there exists a conformal transformation $\varphi:\mathbb{S}^{n}\rightarrow\mathbb{S}^{n}$ such that
$$f=c+\log\det\ud \varphi,\quad c\in \R.$$
		\end{thm A*}

In the range $\gamma>\frac{n}{2}$, $P_{2\gamma}$ can be  also well-defined by (\ref{1.1}),  provided   one sets
\begin{align}\label{1.5}
\frac{1}{\Gamma(l+\frac{n}{2}-\gamma)}=0, \;\; l+\frac{n}{2}-\gamma\in\{0,-1,-2,\cdots\}.
\end{align}
However, in these cases, $ P_{2\gamma} $ is no longer positive definite. Moreover, the exponent $\frac{2n}{n-2\gamma} $ in (\ref{Intro General ineu}) is negative, which necessitates restricting the function $ f $ in (\ref{1.1}) to be positive almost everywhere. In \cite{Frank&Konig&Tang2}, such inequalities are referred to as reverse Sobolev inequalities. An interesting problem is whether the reverse Sobolev inequality still holds. That is, whether the following inequality holds:
			\begin{align}\label{1.6}
			a_{2\gamma}(f)\geq \frac{\Gamma\left(\frac{n}{2}+\gamma\right)}{\Gamma\left(\frac{n}{2}-\gamma\right)}|\S^n|^{\frac{2\gamma}{n}}\|f\|^2_{L^{\frac{2n}{n-2\gamma}}(\S^n)},\;\;
\gamma>\frac{n}{2},\;\; 0<f\in H^{\gamma}(\S^n),
			\end{align}
where we define
	\begin{align}\label{Intro conformal energy}
	 a_{2\gamma}(f): = \int_{\S^n} f P_{2\gamma}f  \ud V_{\S^n}.
	\end{align}

By \eqref{1.1} and (\ref{1.5}), it is straightforward to see that (\ref{1.6}) holds for $\gamma \in \frac{n}{2} + \mathbb{N}$
(here, $\mathbb{N}=\{1,2,3,\cdots\}$), with equality if and only if $f$ is in the linear span of spherical harmonics of degree not exceeding  $\gamma - \frac{n}{2}$ (see \cite{Hang,Frank&Konig&Tang2}).

For other cases, the first reverse Sobolev inequality (for $\gamma = 1$ in $n = 1$) was established by Exner, Harrell, and Loss \cite{Exner&Harrell&Loss} (see also Ai, Chou, and Wei \cite{Ai&Chou&Wei}). For $\gamma = 2$ in $n = 3$, Zhu and Yang \cite{Yang&Zhu} proved the sharp reverse Sobolev inequality using the rearrangement method. In \cite{Hang&Yang1}, Hang and Yang introduced a systematic variational approach and provided a new proof using perturbation techniques (see also the third author \cite{Zhang} for the related Liouville theorem).

We note that the approach by Hang and Yang \cite{Hang&Yang1} relies on the existence of a minimizer, which involves two key technical elements. The first is the approximation of a Sobolev function by functions that vanish near a point (see Lemma 2.3 in \cite{Hang} and Lemma 2.2 in \cite{Hang&Yang1}). The second is the classification of critical exponent equations, as discussed in Li's work ( see \cite{Li}). Using this approach, Hang \cite{Hang} established reverse Sobolev inequalities for odd \(n\) with $\gamma = \frac{n+1}{2}$ or $\frac{n+3}{2}$. Closely following Hang's strategy, Frank, K\"onig, and Tang \cite{Frank&Konig&Tang2} proved the sharp reverse Sobolev inequality for $\gamma \in \left(\frac{n}{2}, \frac{n}{2}+2\right) \setminus \left\{\frac{n}{2} + 1\right\}$. They also proved that the reverse Sobolev inequalities do not hold for $\gamma \in \left(\frac{n}{2}+2, +\infty\right) \setminus \left(\frac{n}{2} + \mathbb{N}\right)$ (see Hang \cite{Hang} for $n$ odd and $\gamma \in \frac{n+5}{2} + \mathbb{N}$).

It is well known that the HLS inequality is equivalent to the Sobolev inequality (\ref{Intro General ineu}). In 2015, Dou and Zhu \cite{Dou&Zhu} posed the question: does this equivalence extend to the reverse HLS inequality and the reverse Sobolev inequality as well? For the reverse HLS inequality in the conformally invariant case, we refer to \cite{Beckner2,Dou&Zhu,Ngo&Nguyen}. For the non-conformally invariant case, see \cite{Carrillo1,Carrillo2}. In 2022, Frank et al. \cite{Frank&Konig&Tang2} said  that there may be no connection between the reverse HLS inequality and the reverse Sobolev inequality, due to the fact that the kernel in the reverse HLS inequality is not positive definite.
However, in Section 3, we demonstrate that the reverse Sobolev inequality for \( \gamma \in \left(\frac{n}{2}, \frac{n}{2}+1\right) \) can indeed be derived from the reversed HLS inequalities, thus partially  answering the question posed by Dou and Zhu \cite{Dou&Zhu}.

For \( \gamma \in \left(\frac{n}{2}+1, \frac{n}{2}+2\right) \), we provide a direct proof of \eqref{1.6} using the vanishing center of mass condition introduced by Hang \cite{Hang}, which is different from that given by Frank and Lieb (see \cite[Lemma B.1]{Frank&Lieb} and \cite[Lemma B.1]{Frank&Lieb2} ) . Our approach does not rely on the strategy of Hang and Yang \cite{Hang&Yang1}, thereby bypassing the two key technical elements they employed.

For readers convenience, we state the  results as the following:
			\begin{thm}\label{Sobolev Thm}
			For  $\gamma\in \left(\frac{n}{2},\frac{n}{2}+1\right)\cup \left(\frac{n}{2}+1,\frac{n}{2}+2\right)$ and $0<f\in H^{\gamma}(\S^n)$,  it holds
			\begin{align}\label{Sobolev Inequality}
		a_{2\gamma}(f) \geq \frac{\Gamma\left(\frac{n}{2}+\gamma\right)}{\Gamma\left(\frac{n}{2}-\gamma\right)}|\S^n|^{\frac{2\gamma}{n}}\|f\|^2_{L^{\frac{2n}{n-2\gamma}}(\S^n)}.
			\end{align}
			The equality holds if and only if here exists a conformal transformation $\varphi: \mathbb{S}^{n}\rightarrow\mathbb{S}^{n}$ such that
$$f=c(\det\ud \varphi)^{\frac{n-2\gamma}{2n}},\; c>0.$$
		\end{thm}
\begin{rem}
 For nonnegative function  $f \in H^{\gamma}(\S^n)$, $\gamma\in \left(\frac{n}{2},\frac{n}{2}+1\right)\cup \left(\frac{n}{2}+1,\frac{n}{2}+2\right)$,  and $\min_{\S^n}f=0$, the Proposition \ref{Positive Pro} states that
 	\begin{align}\label{b1.9}
	\|f\|_{L^{\frac{2n}{n-2\gamma}}(\S^n)}=0.
	\end{align}
By using (\ref{b1.9}) and (\ref{Sobolev Inequality}), one can easy get the inequality (see also \cite{Hang,Frank&Konig&Tang2})
 \begin{align*}
a_{2\gamma}(f)\geq 0.
 \end{align*}
Thus (\ref{Sobolev Inequality}) holds also for nonnegative function.
\end{rem}

On page 6 of \cite{Frank&Konig&Tang2}, Frank et al. posed an open problem for  regarding the search for the quantitative stability of the reverse Sobolev inequality for $\gamma \in \left(\frac{n}{2}, \frac{n}{2}+1\right) \cup \left(\frac{n}{2}+1, \frac{n}{2}+2\right)$.
 We refer the reader   to \cite{Frank2}  For a survey on the subject.
For the quantitative stability of the (usual) Sobolev inequality and  other kinds of functional and geometric inequalities, we refer to \cite{Bianchi&Egnell,Brezis&Lieb, Chen&Kim&Wei, Chen&Frank&Weth,Deng&Sun&Wei,Figalli&Glaudo,Figalli&Neumayer,Figalli&Zhang,Frank,Konig1,Konig2,Lu&Wei,Wei&Wu1,Wei&Wu2}
and the references therein. For explicit lower bounds in the Bianchi-Egnell inequality, we refer to \cite{Chen&Lu&Lu&Tang1,Chen&Lu&Lu&Tang2,Dolbeault&Esteban&Figalli&Frank&Loss}. Since, for \(\gamma > \frac{n}{2}\), the conformal energy \eqref{Intro conformal energy} is not positive definite, the quantitative stability of such inequalities remains unclear. Surprisingly, the method we use to prove the reverse Sobolev inequality for \(\gamma \in \left(\frac{n}{2}+1, \frac{n}{2}+2\right)\) allows us to obtain the quantitative stability of this inequality, which in turn enables us to derive explicit lower bounds. The main result is the following theorem:

		\begin{thm}\label{Stability Thm} Let $\gamma\in \left(\frac{n}{2}+1,\frac{n}{2}+2\right)$ and  $0<f\in H^{\gamma}(\S^n)$.  If $f$ is not a extremal function of reverse Sobolev inequality of order $\gamma$, then  there exist $c>0$ and conformal transformation $\varphi:\mathbb{S}^{n}\rightarrow\mathbb{S}^{n}$ such that
\begin{align*}
a_{2\gamma}\left(f-c\left(\mathrm{\det}\ud \varphi\right)^{\frac{n-2\gamma}{2n}}\right)> 0
\end{align*}
and
			\begin{align*}
			a_{2\gamma}(f)-\frac{\Gamma\left(\frac{n}{2}+\gamma\right)}{\Gamma\left(\frac{n}{2}-\gamma\right)}|\S^n|^{\frac{2\gamma}{n}}\|f\|^2_{L^{\frac{2n}{n-2\gamma}}(\S^n)}\geq a_{2\gamma}\left(f-c\left(\mathrm{\det}\ud \varphi\right)^{\frac{n-2\gamma}{2n}}\right).
			\end{align*}
		\end{thm}

Finally, we shall use these Sobolev inequalities  along with conformally covariant boundary operators to establish Sobolev trace inequalities on the unit ball $\mathbb{B}^{n+1} = \{x \in \mathbb{R}^{n+1}: |x| < 1\}$.
We note that the first inequality of this type was established by Escobar \cite{es1} (see also Beckner \cite{Beckner}). It states that for $f \in C^{\infty}(\mathbb{S}^{n})$ and $n \geq 2$,
  \begin{equation}\label{1.9}
  \begin{split}
\frac{n-1}{2}|\S^n|^{\frac{1}{n}}\|f\|^2_{L^{\frac{2n}{n-1}}(\S^n)} \leq
\int_{\mathbb{B}^{n+1}}|\nabla u|^{2}dx+\frac{n-1}{2}\|f\|^2_{L^{2}(\S^n)},
\end{split}
\end{equation}
where $u$ is  a smooth extension of $f$ to $\mathbb{B}^{n+1}$. In  the limiting case $n=1$, one has the    Lebedev-Milin  inequality (\cite{Lebedev})
  \begin{equation}\label{1.10}
  \begin{split}
\log \fint_{\mathbb{S}^{1}}e^{f}\ud V_{{\S^1}}\leq
\frac{1}{4\pi}
\int_{\mathbb{B}^{2}}|\nabla u|^{2}dx+\fint_{\mathbb{S}^{1}}f\ud V_{{\S^1}}.
\end{split}
\end{equation}
Such inequalities have been generalized to higher-order derivatives on \(\mathbb{B}^{n+1}\) by several authors (see \cite{Ache&Chang,Chen&Zhang,Ngo&Nguyen&Pham,yang}
) and to the CR setting by Frank et al. \cite{Frank2} and Flynn, Lu, and the second author \cite{Flynn&Lu&Yang1}.

In \cite{Case2}, Case introduced a series of boundary operators on half-spaces and used them to establish a family of sharp Sobolev trace inequalities on half-spaces. However, it remains unclear whether the boundary operators in \cite{Case2} satisfy the conformal covariance property.  Additionally, it seems that the method of introducing boundary operators on half-spaces in \cite{Case2} does not apply to the case of \(\mathbb{B}^{n+1}\). Recently, Flynn, Lu, and the second author \cite{Flynn&Lu&Yang1,Flynn&Lu&Yang} introduced another type of  boundary operators for Poincar\'e-Einstein manifolds and the complex ball, which satisfy the conformal covariance property. Using these boundary operators and  scattering theory (see \cite{ma,Graham&Zworski}), they also established higher-order (usual) Sobolev trace inequalities on Poincar\'e-Einstein manifolds and the complex ball. For other  conformal boundary operators associated with the GJMS operators on manifolds, we refer to \cite{Case1,Case&Luo}.
In this paper, we will demonstrate that the approach,  introduced in  \cite{Flynn&Lu&Yang1,Flynn&Lu&Yang}, together with reverse Sobolev inequalities, can be also applicable to obtain Sobolev trace inequalities on \(\mathbb{B}^{n+1}\) for higher-order derivatives.

To state our results, let us agree to some conventions.
 Throughout the paper, we adopt the following notation and notions:
\begin{itemize}
	\item We consider  $\gamma\in (0,+\infty)\backslash \mathbb{N}$, where $\mathbb{N}$ denotes the set of  positive integers, specifically given by
	\begin{align*}
	\mathbb{N}=\{1,2,3,\cdots\}.
	\end{align*}
 Let $\lfloor \gamma\rfloor$ denote the integer part of $\gamma$ and  $[\gamma]=\gamma-\lfloor \gamma\rfloor$  the fractional part. Additionally,  we often work with the constant
	\begin{align}\label{Intro c gamma}
	c_{\gamma}=\frac{2^{2\gamma}\Gamma(\gamma)}{\Gamma(-\gamma)}.
	\end{align}
	 \item  Let $ g_{\B}$ be the Poincar\'e metric on  $\mathbb{B}^{n+1}$:
	\begin{align*}
	g_{\B}=\frac{4|\ud x|^2}{(1-|x|^2)^2}=\frac{4(dr^{2}+r^{2}d\theta^{2})}{(1-r^{2})^{2}},
	\end{align*}
	where $(r,\theta)$ is the polar coordinate  on $\mathbb{B}^{n+1}$ and  $d\theta^{2}$ is the canonical
	spherical metric.
	\item  From the classical scattering theory, we need to introduce the two type defining functions as follows: the geodesic normal defining function
	\begin{align*}
	\rho=\frac{2(1-r)}{1+r},
	\end{align*}
	and the flat defining function
	\begin{align*}
	\rho_0=\frac{1-r^2}{2}.
	\end{align*}
	\item Denote by $\Delta_{+}$  the Laplace-Beltrami operator on $(\mathbb{B}^{n+1}, g_{\B})$. One has
	\begin{align}\label{laplace poin oprtaor}
	\Delta_{+}=\frac{1-|x|^2}{4}\left[(1-|x|^2)\Delta_{\R^{n+1}}+2(n-1)\sum_{i=1}^{n+1}x_i\partial_{x_i}\right].
	\end{align}
	For $s=\frac{n}{2}+\gamma$, we let
	\begin{align*}
	D_s=-\Delta_{+}-s(n-s),\qquad L_{2k}^{+}=\prod_{j=0}^{k-1}D_{s-2j}, \qquad k=\lfloor\gamma\rfloor+1
	\end{align*}
	and
	\begin{align*}
	L_{2k}=\rho^{-\frac{n}{2}+\gamma-2k}\circ L^{+}_{2k}\circ \rho^{\frac{n}{2}-\gamma}.
	\end{align*}
	\item For $\gamma\in (0,+\infty)$, we consider the following function space
	\begin{align*}
	\mathcal{C}^{2\gamma}(\overline{\B^{n+1}})=C^{\infty}_{\mathrm{even}}(\overline{\B^{n+1}})+\rho^{2[\gamma]}C^{\infty}_{\mathrm{even}}(\overline{\B^{n+1}}),
	\end{align*}
	where $f\in C^{\infty}_{\mathrm{even}}(\overline{\B^{n+1}})$ indicates that $f\in C^{\infty}(\overline{\B^{n+1}})$ and possess a Taylor expansion of the form
	\begin{align*}
	f=f_0+\rho^2f_2+\rho^4 f_4+\cdots,\quad f_i\in C^{\infty}(\S^n)\quad \mathrm{for}\quad i=0,2,\cdots.
	\end{align*}
	Especially, it is noteworthy that if $[\gamma]=\frac{1}{2}$,then  $\mathcal{C}^{2\gamma}(\overline{\B^{n+1}})=C^{\infty}(\overline{\B^{n+1}})$.
	\item We adopt the notation $O(\rho^{\alpha})$ to indicate that  there exist smooth functions $g_0,g_1,g_2,\cdots\in C^{\infty}(\S^n)$ such that
	\begin{align*}
	O(\rho^{\alpha})=\rho^{\alpha}(g_0+\rho^{2}g_1+\rho^{4}g_2\cdots).
	\end{align*}
	\item A defining function $\hat{\rho}$ of $\mathbb{B}^{n+1}$ is called $\gamma$-admissible provided
	\begin{equation}
	\frac{\rho} {\hat{\rho}}= 1 + \sum_{j=1}^{\lfloor \gamma \rfloor} \rho_{(2j)} {\hat{\rho}}^{2j} + \Phi {\hat{\rho}}^{2 \gamma} + o({\hat{\rho}}^{2\gamma}) \text{ as } {\hat{\rho}}\to0,
	\label{eq:gamma-admissible-defining-function}
	\end{equation}
	where $\rho_{(2j)}(1\leq j\leq \lfloor \gamma \rfloor),\Phi \in C^{\infty}(\mathbb{S}^{n})$.
	Clearly, if $\hat{\rho}$ is $\gamma$-admissible, then
	\begin{equation}
	\frac{{\hat{\rho}}}{\rho} = 1 + \sum_{j=1}^{\lfloor \gamma \rfloor} \hat\rho_{(2j)} \rho^{2j} + \hat\Phi \rho^{2 \gamma} + o(\rho^{2\gamma}) \text{ as } \rho \to0,
	\end{equation}
	for some $ \hat\rho_{(2j)}(1\leq j\leq \lfloor \gamma \rfloor),\hat\Phi \in C^{\infty}(\mathbb{S}^{n})$.
\end{itemize}

With the above notation, we begin by introducing the boundary operators in a simple case, namely when the index is small. For a given $ U \in \mathcal{C}^{2\gamma}(\overline{\B^{n+1}}) $, we define the boundary operators associated with the metric $ g = \rho^2 g_{\B} $ for small indices, following the approach in \cite{Flynn&Lu&Yang1,Flynn&Lu&Yang}, as follows:

\begin{itemize}
\item $ B_{ 0}^{2\gamma}(U)=U|_{\rho=0}$;
  \item  for $1\leq j\leq \lfloor\gamma/2\rfloor$, define
  \begin{align*}
B^{2\gamma}_{2j}(U)&=\frac{1}{b_{2j}}\rho^{-\frac{n}{2}+\gamma-2j}\prod_{\ell=0}^{j-1}D_{s-2l} \prod_{\ell=\lfloor\gamma\rfloor-j+1}^{\lfloor\gamma\rfloor}D_{s-2l} \left(\rho^{\frac{n}{2}-\gamma}U\right)\Big|_{\rho=0};
  \end{align*}
  \item for $0\leq j\leq \lfloor\gamma\rfloor-\lfloor\gamma/2\rfloor-1$, define
  \begin{align*}
B^{2\gamma}_{2j+2[\gamma]}(U)&=\frac{-1}{b_{2j+2[\gamma]}}\rho^{-\frac{n}{2}+\gamma-2j-2[\gamma]} \prod_{l=0}^{j}D_{s-2l} \prod_{l=\lfloor\gamma\rfloor-j+1}^{\lfloor\gamma\rfloor}D_{s-2l}\left(\rho^{\frac{n}{2}-\gamma}U\right)\Big|_{\rho=0},
\end{align*}
\end{itemize}
where
	\begin{align*}
b_{2j}&=4^{2j}j!\frac{\Gamma(j+1-[\gamma])}{\Gamma(1-[\gamma])}\frac{\Gamma(\gamma+1-j)}{\Gamma(\gamma+1-2j)}\frac{\Gamma(\lfloor \gamma\rfloor+1-j)}{\Gamma(\lfloor \gamma\rfloor+1-2j)};\\
b_{2j+2[\gamma]}&=-4^{2j}j!\frac{\Gamma(j+1+[\gamma])}{\Gamma([\gamma])}\frac{\Gamma(\lfloor \gamma\rfloor+1-j)}{\Gamma(\lfloor \gamma\rfloor-2j)}\frac{\Gamma(\lfloor \gamma\rfloor+1-j-[\gamma])}{\Gamma(\lfloor \gamma\rfloor+1-2j-[\gamma])}
\end{align*}
are  chosen such that
\begin{align*}
B^{2\gamma}_{2j}(\rho^{2j})=B^{2\gamma}_{2j+2[\gamma]}(\rho^{2j+2[\gamma]})=1.
\end{align*}
The definition of  boundary operators with large indices is  in terms of the solution $L_{2k} U = 0$.
 This is done in subsection \ref{sec large index}.

 \begin{convention}
		Let  $g=\rho^2g_{\B}$. In this framework, our boundary operators are defined as the differential operators associated with the metric $g$, namely,
		\begin{align*}
		B^{2\gamma}_{2j}:=B^{2\gamma,g}_{2j},\qquad B^{2\gamma}_{2j+2[\gamma]}:=B^{2\gamma,g}_{2j+2[\gamma]},\qquad P_{2\gamma}:=P^{g|_{\S^n}}_{2\gamma},
		\end{align*}
		where $P_{2\gamma}$ denotes the GJMS operator on the boundary. For clarity and simplicity, we often omit the superscript
		 $g$ in our notation.
		\end{convention}

 Let $\widehat{\rho}=e^{\tau}\rho$ be another $\gamma$-admissible defining function.   Denote by  $\widehat{B}_{2j}^{2\gamma}$ and $\widehat{B}^{2\gamma}_{2j+2[\gamma]}$  the boundary operators associated with $(\mathbb{B}^{n+1}, \widehat{g}=e^{2\tau}g)$. We have the following theorem:
		\begin{thm}\label{Conformal Thm}
			Let $\gamma\in (0,+\infty)\backslash\mathbb{N}$ and $U\in 	\mathcal{C}^{2\gamma}(\overline{\B^{n+1}})$.  For $0\leq j\leq \lfloor \gamma\rfloor$,   we have \begin{align}\label{Conformal Thm property 1}
				\widehat{B}^{2\gamma}_{2j}(U)=&e^{\left(-\frac{n}{2}+\gamma-2j\right)\tau|_{\S^n}}B^{2\gamma}_{2j}\left(e^{\left(\frac{n}{2}-\gamma\right)\tau}U\right);\\
\label{Conformal Thm property 2}				\widehat{B}^{2\gamma}_{2j+2[\gamma]}(U)=&e^{\left(-\frac{n}{2}+\gamma-2j-2[\gamma]\right)\tau|_{\S^n}}B^{2\gamma}_{2j+2[\gamma]}
\left(e^{\left(\frac{n}{2}-\gamma\right)\tau}U\right).
			\end{align}
		\end{thm}

By the classical Caffarelli-Silvestre theory in \cite{Caffarelli&Silvestre}, it is well-known that the boundary operator recovers the fractional GJMS operator on the boundary. The following theorem demonstrates that our boundary operators exhibit similar properties.

		\begin{thm}\label{Intrinsic Thm}
			Let $\gamma\in (0,+\infty)\backslash\mathbb{N}$ and $k=\lfloor \gamma\rfloor+1$. For any $U\in 	\mathcal{C}^{2\gamma}(\overline{\B^{n+1}})$ satisfying $L_{2k}U=0$, it holds
			\begin{align*}
			B^{2\gamma}_{2\gamma-2j}(U)=&\frac{1}{c_{\gamma-2j}}P_{{2\gamma-4j}}B^{2\gamma}_{2j}(U) \qquad\qquad\qquad~~\qquad\mathrm{for}\qquad 0\leq j\leq \lfloor \gamma/2\rfloor;\\
					B^{2\gamma}_{2\lfloor \gamma\rfloor-2j}(U)=&\frac{1}{c_{\lfloor \gamma\rfloor-[\gamma]-2j}}P_{2\lfloor \gamma\rfloor-2[\gamma]-4j}B^{2\gamma}_{2[\gamma]+2j}(U) \qquad\mathrm{for}\qquad 0\leq j\leq \lfloor \gamma\rfloor-\lfloor \gamma/2\rfloor-1,
			\end{align*}
where $c_{\gamma-2j}, c_{\lfloor \gamma\rfloor-[\gamma]-2j}$ are given by \eqref{Intro c gamma}.
		\end{thm}
		
		Next, we consider  the associated Dirichlet form
		\begin{align*}
		\mathcal{Q}_{2\gamma}(U,V)=&\int_{\B^{n+1}}\rho^{1-2[\gamma]}UL_{2k}V\ud V_{g}-\sum_{j=0}^{\lfloor \gamma/2\rfloor}\sigma_{j,\gamma}\int_{\S^n}B^{2\gamma}_{2j}UB^{2\gamma}_{2\gamma-2j}V\ud V_{\S^n}\\
		-&\sum_{j=\lfloor \gamma/2\rfloor+1}^{\lfloor \gamma\rfloor}\sigma_{j,\gamma}\int_{\S^n}B^{2\gamma}_{2j}VB^{2\gamma}_{2\gamma-2j}U\ud V_{\S^n},
		\end{align*}
	where
	\begin{align*}
		\sigma_{j,\gamma}=\begin{cases}
		\displaystyle 2^{2\lfloor \gamma\rfloor+1}j!(\lfloor \gamma\rfloor-j)!\frac{\Gamma(\gamma+1-j)}{\Gamma(\gamma-2j)}\frac{\Gamma(j+1-[\gamma])}{\Gamma(2j+1-\gamma)}, \quad 0\leq j\leq \lfloor \gamma/2\rfloor;\\
		\displaystyle -2^{2\lfloor \gamma\rfloor+1}j!(\lfloor \gamma\rfloor-j)!\frac{\Gamma(\gamma+1-j)}{\Gamma(\gamma-2j)}\frac{\Gamma(j+1-[\gamma])}{\Gamma(2j+1-\gamma)},\quad \lfloor \gamma/2\rfloor+1\leq j\leq \lfloor \gamma\rfloor.
		\end{cases}
	\end{align*}
	
	Using Theorem \ref{Conformal Thm}, we establish that the Dirichlet form $ \mathcal{Q}_{2\gamma}(U, V) $ is conformally invariant in the following manner: if $ \hat{\rho} = e^{\tau} \rho $, then we have
	\begin{align*}
		\hat{\mathcal{Q}}_{2\gamma}(U,V)=\mathcal{Q}_{2\gamma}(e^{\left(\frac{n}{2}-\gamma\right)\tau}U,e^{\left(\frac{n}{2}-\gamma\right)\tau}V).
	\end{align*}
Another significant property of $ \mathcal{Q}_{2\gamma}(U, V) $ is its symmetry, which is deeply rooted in the construction of the conformal boundary operators. The proof of this symmetry is intricate and highly nontrivial.

		\begin{thm}\label{Symmetric Thm}
			Let $U, V\in 	\mathcal{C}^{2\gamma}(\overline{\B^{n+1}})$, there holds
			\begin{align*}
				\mathcal{Q}_{2\gamma}(U,V)=\mathcal{Q}_{2\gamma}(V,U)
			\end{align*}
		\end{thm}
		
		 Using the results mentioned above,  we can derive the sharp Sobolev trace inequality. With this, the Dirichlet form $ \mathcal{Q}_{2\gamma}(U, U) $can be interpreted as the conformal energy.
		
		\begin{thm}\label{Functional inequ Thm}
			Let $\gamma\in (0,+\infty)\backslash\mathbb{N}$, $k=\lfloor \gamma\rfloor+1$ and $U\in 	\mathcal{C}^{2\gamma}(\overline{\B^{n+1}})$. It holds
			\begin{align*}
					\mathcal{Q}_{2\gamma}(U,U)\geq &
					\sum_{j=0}^{\lfloor \gamma/2\rfloor}\upzeta_{j,\gamma}\int_{\S^n}B^{2\gamma}_{2j}UP_{2\gamma-4j}B^{2\gamma}_{2j}U\ud V_{\S^n}\\
					+&\sum_{j=\lfloor \gamma/2\rfloor+1}^{\lfloor \gamma\rfloor}\upzeta_{j,\gamma}\int_{\S^n}B^{2\gamma}_{2\gamma-2j}UP_{4j-2\gamma}B^{2\gamma}_{2\gamma-2j}U\ud V_{\S^n},
			\end{align*}
			where
			\begin{align*}
					\upzeta_{j,\gamma}=\begin{cases}
				\displaystyle 2^{4j-2[\gamma]+1}j!(\lfloor \gamma\rfloor-j)!\frac{\Gamma(\gamma+1-j)}{\Gamma(\gamma+1-2j)}\frac{\Gamma(j+1-[\gamma])}{\Gamma(\gamma-2j)}, \quad 0\leq j\leq \lfloor \gamma/2\rfloor;\\
				\displaystyle 2^{4\gamma-2[\gamma]+1}j!(\lfloor \gamma\rfloor-j)!\frac{\Gamma(\gamma+1-j)}{\Gamma(2j-\gamma)}\frac{\Gamma(j+1-[\gamma])}{\Gamma(2j+1-\gamma)},\quad \lfloor \gamma/2\rfloor+1\leq j\leq \lfloor \gamma\rfloor.
				\end{cases}
			\end{align*}
			Moreover, the equality holds if and only if $L_{2k}(U)=0$.
		\end{thm}

By employing Theorem \ref{Functional inequ Thm} along with the (usual and reverse) Sobolev inequalities, we can provide a comprehensive version of the Sobolev trace inequality for all \( \gamma \in (0, +\infty) \setminus \mathbb{N} \). To ensure clarity, we introduce some new notations. For $\gamma\in \left(0,+\infty\right)\backslash\mathbb{N}$	, we introduce two integers $j_1=\lfloor \frac{\gamma}{2}-\frac{n}{4}\rfloor$ and $j_{2}=\lfloor \frac{\gamma}{2}+\frac{n}{4}\rfloor+1$.	If $\gamma\in \left(\frac{n}{2},+\infty\right)\backslash\mathbb{N}$, we define the function space
		\begin{align*}
			\mathscr{N}^{2\gamma}=&\left\{U\in 	\mathcal{C}^{2\gamma}(\overline{\B^{n+1}})\Big|B^{2\gamma}_{2j}U=0, ~\mathrm{for}~  j<j_1 ,  ~ \mathrm{and}~ B^{2\gamma}_{2j_1}U> 0\right\}\cap\\
			&\left\{U\in 	\mathcal{C}^{2\gamma}(\overline{\B^{n+1}})\Big|B^{2\gamma}_{2\gamma-2j}U=0 , ~\mathrm{for}~ j_2<j\leq \lfloor \gamma\rfloor, ~ \mathrm{and}~ B^{2\gamma}_{2\gamma-2j_2}U>0\right\}.
		\end{align*}
	if $\gamma\in \left(0, \frac{n}{2}\right)\backslash\mathbb{N}$, then we can define
	\begin{align*}
			\mathscr{N}^{2\gamma}=\mathcal{C}^{2\gamma}(\overline{\B^{n+1}}).
	\end{align*}
		
To this end, we have the following corollary:		
		\begin{cor}\label{Sobolev Trace Thm}
		Suppose $\gamma\in (0,+\infty)\backslash\mathbb{N}$ and let $k=\lfloor \gamma\rfloor+1$. We have the following inequalities:
		\begin{itemize}	
			\item[(I)]  $[\gamma]=\frac{1}{2}$ and $n$ is odd.
			\begin{itemize}
				\item[(1)] If $\gamma-\frac{n}{2}\in 2\mathbb{N}\cup\{0\}$, then  for any $U\in 	C^{\infty}(\overline{\B^{n+1}})$, it holds
				\begin{align*}
				\mathcal{Q}_{2\gamma}(U,U)\geq& 2n!|\S^n|\upzeta_{j_1,\gamma}\log \left(\fint_{\S^n} e^{B^{2\gamma}_{2j_1}(U)-\overline{B^{2\gamma}_{2j_1}(U)}}  \ud V_{{\S^n}}\right)\\
				+&\sum_{j=j_1+1}^{\lfloor \gamma/2\rfloor}\frac{\Gamma\left(\frac{n}{2}+\gamma-2j\right)}{\Gamma\left(\frac{n}{2}-\gamma+2j\right)}|\S^n|^{\frac{2\gamma-4j}{n}}
				\upzeta_{j,\gamma}\|B^{2\gamma}_{2j}U\|^2_{L^{\frac{2n}{n-2\gamma+4j}}(\S^n)}\\
				+&\sum_{j=\lfloor \gamma/2\rfloor+1}^{j_2-1}\frac{\Gamma\left(\frac{n}{2}-\gamma+2j\right)}{\Gamma\left(\frac{n}{2}+\gamma-2j\right)}|\S^n|^{\frac{4j-2\gamma}{n}}
				\upzeta_{j,\gamma}\|B^{2\gamma}_{2\gamma-2j}U\|^2_{L^{\frac{2n}{n+2\gamma-4j}}(\S^n)}.
				\end{align*}
				The equality holds if and only if $U$ is the unique solution of
				\begin{align*}
				\begin{cases}
				\displaystyle L_{2k}U=0 \qquad&\mathrm{in}\qquad \B^{n+1},\\
				\displaystyle B^{2\gamma}_{2j}(U)\in \oplus_{l=0}^{\lfloor \gamma-2j-\frac{n}{2}\rfloor}\mathscr{H}_{l} \qquad&\mathrm{for}\qquad 0\leq j<j_1,\\
				\displaystyle B^{2\gamma}_{2j_1}(U)=a_{2j_1}+\log \det\ud \varphi_{2j_1} \qquad&\mathrm{for}\qquad  j=j_1,\\
				\displaystyle B^{2\gamma}_{2j}(U)=a_{2j}(\det\ud \varphi_{2j})^{\frac{n-2\gamma+4j}{2n}} \qquad&\mathrm{for}\qquad  j_1<j\leq
				\lfloor \gamma/2\rfloor,\\
				\displaystyle B^{2\gamma}_{2\gamma-2j}(U)=a_{2\gamma-2j}(\det\ud \varphi_{2\gamma-2j})^{\frac{n+2\gamma-4j}{2n}} \qquad&\mathrm{for}\qquad  \lfloor \gamma/2\rfloor+1\leq j<j_2,\\
				\displaystyle B^{2\gamma}_{2\gamma-2j}(U)\in\oplus_{l=0}^{\lfloor 2j-\gamma-\frac{n}{2}\rfloor}\mathscr{H}_{l} \qquad&\mathrm{for}\qquad  j_2\leq j\leq\lfloor \gamma\rfloor.
				\end{cases}
				\end{align*}
				\item[(2)] If $\gamma+\frac{n}{2}\in 2\mathbb{N}$, then  for any $U\in 	C^{\infty}(\overline{\B^{n+1}})$, it holds
				\begin{align*}
				\mathcal{Q}_{2\gamma}(U,U)\geq& 2n!|\S^n|\chi_{\{j_2-1,\lfloor \gamma\rfloor\}}\upzeta_{j_2-1,\gamma}\log \left(\fint_{\S^n} e^{B^{2\gamma}_{2\gamma-2j_2+2}(U)-\overline{B^{2\gamma}_{2\gamma-2j_2+2}(U)}}  \ud V_{{\S^n}}\right)\\
				+&\sum_{j=\max\{j_1,0\}+1}^{\lfloor \gamma/2\rfloor}\frac{\Gamma\left(\frac{n}{2}+\gamma-2j\right)}{\Gamma\left(\frac{n}{2}-\gamma+2j\right)}|\S^n|^{\frac{2\gamma-4j}{n}}
				\upzeta_{j,\gamma}\|B^{2\gamma}_{2j}U\|^2_{L^{\frac{2n}{n-2\gamma+4j}}(\S^n)}\\
				+&\sum_{j=\lfloor \gamma/2\rfloor+1}^{\min\{j_2-2,\lfloor \gamma\rfloor\}}\frac{\Gamma\left(\frac{n}{2}-\gamma+2j\right)}{\Gamma\left(\frac{n}{2}+\gamma-2j\right)}|\S^n|^{\frac{4j-2\gamma}{n}}
				\upzeta_{j,\gamma}\|B^{2\gamma}_{2\gamma-2j}U\|^2_{L^{\frac{2n}{n+2\gamma-4j}}(\S^n)},
				\end{align*}
where $\chi_{\{j_2-1,\lfloor \gamma\rfloor\}}=1$ if $j_2-1\leq\lfloor \gamma\rfloor$, otherwise $\chi_{\{j_2-1,\lfloor \gamma\rfloor\}}=0$.
				The equality holds if and only if $U$ is the unique solution of
				\begin{align*}
				\begin{cases}
				\displaystyle L_{2k}U=0 \quad&\mathrm{in}\quad \B^{n+1},\\
				\displaystyle B^{2\gamma}_{2j}(U)\in \oplus_{l=0}^{\lfloor \gamma-2j-\frac{n}{2}\rfloor}\mathscr{H}_{l} \quad&\mathrm{for}\quad 0\leq j\leq j_1,\\
				\displaystyle B^{2\gamma}_{2j}(U)=a_{2j}(\det\ud \varphi_{2j})^{\frac{n-2\gamma+4j}{2n}} \quad&\mathrm{for}\quad  \max\{j_1+1,0\}\leq j\leq
				\lfloor \gamma/2\rfloor,\\
				\displaystyle B^{2\gamma}_{2\gamma-2j}(U)=a_{2\gamma-2j}(\det\ud \varphi_{2\gamma-2j})^{\frac{n+2\gamma-4j}{2n}}\quad&\mathrm{for}\quad  \lfloor \gamma/2\rfloor+1\leq j\leq\min\{j_2-2,\lfloor \gamma\rfloor\},\\
				\displaystyle B^{2\gamma}_{2\gamma-2j_2+2}(U)=a_{2\gamma-2j_2+2}+\log\det\ud \varphi_{2\gamma-2j_2+2}\quad&\mathrm{for}\quad   j=j_2-1\;\;
\mathrm{when}\;\;j_2-1\leq \lfloor \gamma\rfloor,\\
				\displaystyle B^{2\gamma}_{2\gamma-2j}(U)\in\oplus_{l=0}^{\lfloor 2j-\gamma-\frac{n}{2}\rfloor}\mathscr{H}_{l} \quad&\mathrm{for}\quad  j_2\leq j\leq\lfloor \gamma\rfloor.
				\end{cases}
				\end{align*}
				
				\item[(3)]  If $\gamma-\frac{n}{2}\not\in 2\mathbb{N}\cup\{0\}$ and $\gamma+\frac{n}{2}\not\in 2\mathbb{N}$ ,then  for any $U\in 	C^{\infty}(\overline{\B^{n+1}})$, it holds
				\begin{align*}
				\mathcal{Q}_{2\gamma}(U,U)\geq
				&\sum_{j=\max\{j_1,0\}+1}^{\lfloor \gamma/2\rfloor}\frac{\Gamma\left(\frac{n}{2}+\gamma-2j\right)}{\Gamma\left(\frac{n}{2}-\gamma+2j\right)}|\S^n|^{\frac{2\gamma-4j}{n}}
				\upzeta_{j,\gamma}\|B^{2\gamma}_{2j}U\|^2_{L^{\frac{2n}{n-2\gamma+4j}}(\S^n)}\\
				+&\sum_{j=\lfloor \gamma/2\rfloor+1}^{\min\{j_2-1,\lfloor \gamma\rfloor\}}\frac{\Gamma\left(\frac{n}{2}-\gamma+2j\right)}{\Gamma\left(\frac{n}{2}+\gamma-2j\right)}|\S^n|^{\frac{4j-2\gamma}{n}}
				\upzeta_{j,\gamma}\|B^{2\gamma}_{2\gamma-2j}U\|^2_{L^{\frac{2n}{n+2\gamma-4j}}(\S^n)}.
				\end{align*}
				The equality holds if and only if $U$ is the unique solution of
				\begin{align*}
				\begin{cases}
				\displaystyle L_{2k}U=0 \qquad&\mathrm{in}\qquad \B^{n+1},\\
				\displaystyle B^{2\gamma}_{2j}(U)\in \oplus_{l=0}^{\lfloor \gamma-2j-\frac{n}{2}\rfloor}\mathscr{H}_{l},\qquad&\mathrm{for}\qquad 0\leq j\leq  j_1,\\
				\displaystyle B^{2\gamma}_{2j}(U)=a_{2j}(\det\ud \varphi_{2j})^{\frac{n-2\gamma+4j}{2n}},\qquad&\mathrm{for}\qquad   \max\{j_1+1,0\}\leq j\leq
				\lfloor \gamma/2\rfloor\\
				\displaystyle B^{2\gamma}_{2\gamma-2j}(U)=a_{2\gamma-2j}(\det\ud \varphi_{2\gamma-2j})^{\frac{n+2\gamma-4j}{2n}}\qquad&\mathrm{for}\qquad  \lfloor \gamma/2\rfloor+1\leq j\leq \min\{j_2-1,\lfloor \gamma\rfloor\}\\
				\displaystyle B^{2\gamma}_{2\gamma-2j}(U)\in\oplus_{l=0}^{\lfloor 2j-\gamma-\frac{n}{2}\rfloor}\mathscr{H}_{l}.\qquad&\mathrm{for}\qquad  j_2\leq j\leq\lfloor \gamma\rfloor.
				\end{cases}
				\end{align*}
			\end{itemize}
			\item[(II)]    $[\gamma]\not=\frac{1}{2}$ or $n$ is even. For any  $U\in 	\mathscr{N}^{2\gamma}$,  it holds
			\begin{align*}
			\mathcal{Q}_{2\gamma}(U,U)\geq & \sum_{j=\max\{j_1,0\}}^{\lfloor \gamma/2\rfloor}\frac{\Gamma\left(\frac{n}{2}+\gamma-2j\right)}{\Gamma\left(\frac{n}{2}-\gamma+2j\right)}|\S^n|^{\frac{2\gamma-4j}{n}}
			\upzeta_{j,\gamma}\|B^{2\gamma}_{2j}U\|^2_{L^{\frac{2n}{n-2\gamma+4j}}(\S^n)}\\
			+&\sum_{j=\lfloor \gamma/2\rfloor+1}^{\min\{j_2, \lfloor \gamma\rfloor\}}\frac{\Gamma\left(\frac{n}{2}-\gamma+2j\right)}{\Gamma\left(\frac{n}{2}+\gamma-2j\right)}|\S^n|^{\frac{4j-2\gamma}{n}}
			\upzeta_{j,\gamma}\|B^{2\gamma}_{2\gamma-2j}U\|^2_{L^{\frac{2n}{n+2\gamma-4j}}(\S^n)}.
			\end{align*}
			Moreover, the equality holds if and only if $U$ is the unique solution to
			\begin{align*}
			\begin{cases}
			\displaystyle L_{2k}U=0 \qquad&\mathrm{in}\qquad \B^{n+1},\\
			\displaystyle B^{2\gamma}_{2j}(U)=0\qquad&\mathrm{for}\qquad 0\leq j<\max\{j_1,0\},\\
			\displaystyle B^{2\gamma}_{2j}(U)=a_{2j}(\det\ud \varphi_{2j})^{\frac{n-2\gamma+4j}{2n}}\qquad&\mathrm{for}\qquad \max\{j_1,0\}\leq j\leq \lfloor \gamma/2\rfloor,\\
			\displaystyle
			B^{2\gamma}_{2\gamma-2j}(U)=a_{2\gamma-2j}(\det\ud \varphi_{2\gamma-2j})^{\frac{n+2\gamma-4j}{2n}}\qquad&\mathrm{for}\qquad \lfloor \gamma/2\rfloor+1\leq  j\leq \min\{j_2, \lfloor \gamma\rfloor\},\\
			\displaystyle B^{2\gamma}_{2\gamma-2j}(U)=0\qquad&\mathrm{for}\qquad  \min\{j_2, \lfloor \gamma\rfloor\}< j\leq \lfloor \gamma\rfloor.
			\end{cases}
			\end{align*}
			where $a_{2j}, a_{2\gamma-2j}\in \R_{+}$ and $\varphi_{2j}, \varphi_{2\gamma-2j}$ are two conformal transformations on $\S^n$.

		\end{itemize}
We must point out that in all cases, if no value of $j$ satisfies the above constraint (e.g., if $\gamma < n/2$, then there is no $j$ such that $ \min\{j_2, \lfloor \gamma \rfloor\} < j \leq \lfloor \gamma \rfloor$), then the corresponding boundary condition does not exist when equality holds.
	\end{cor}
		
		The paper is organized as follows. In Section \ref{Sec 2}, we introduce some foundational material on the hypergeometric function and the Funk-Hecke formula, which will be used in scattering theory and in deriving the formula for $P^{-1}_{2\gamma}$. In Section \ref{Sec 3}, we prove Theorem \ref{Sobolev Thm} by considering two cases: one for $\gamma \in \left(\frac{n}{2}, \frac{n}{2}+1\right)$ and the other for $\gamma \in \left(\frac{n}{2}+1, \frac{n}{2}+2\right) $. We also provide a proof of the  quantitative stability of reverse Sobolev inequalities of order  $\gamma \in \left(\frac{n}{2}+1, \frac{n}{2}+2\right)$, i.e. Theorem \ref{Stability Thm}. In Section \ref{Sec 4}, we first present explicit formulas for the boundary operator with small indices, followed by establishing the Dirichlet solvability theorem. Based on this theorem, we then define the boundary operator for larger indices  and  prove Theorem \ref{Conformal Thm}, \ref{Intrinsic Thm}, \ref{Symmetric Thm} and \ref{Functional inequ Thm}. Finally, combining all these results, we conclude with a complete Sobolev trace inequality.

			\section{Preliminaries}\label{Sec 2}
			\subsection{Hypergeometric function}

		First, we recall some basic properties of the hypergeometric function. For further details on these functions, we refer the reader to \cite{Gradshteyn&Ryzhik}.
 We use the notation $ F(a,b;c;z) $ to represent
			\begin{equation}\label{Pre 1.1}
			\begin{split}
			F(a,b;c;z)=\sum^{\infty}_{k=0}\frac{(a)_{k}(b)_{k}}{(c)_{k}}\frac{z^{k}}{k!},
			\end{split}
			\end{equation}
			where $ c \neq 0, -1, -2, \dots $ and $ (a)_k $ denotes the rising Pochhammer symbol, defined by	
			$$
			(a)_{0}=1,\;(a)_{k}=a(a+1)\cdots(a+k-1), \;k\geq1.
			$$
			If either $a$ or $b$ is a nonpositive integer, then the series terminates, and the function reduces to a polynomial.
			
		Next, we list only the properties of the hypergeometric function that will be used in the remainder of the paper.
			\begin{itemize}
				\item The hypergeometric function $F(a,b;c;z)$ satisfies  the hypergeometric differential equation
				\begin{equation}\label{Pre 1.2}
				\begin{split}
				z(1-z)F''+(c-(a+b+1)z)F'-abF=0.
				\end{split}
				\end{equation}
				
				\item If $\textrm{Re} (c-a-b)>0$, then $F(a,b;c;1)$ exists and
				\begin{equation}\label{Pre 1.3}
				\begin{split}
				F(a,b;c;1)=\frac{\Gamma(c)\Gamma(c-a-b)}{\Gamma(c-a)\Gamma(c-b)}.
				\end{split}
				\end{equation}

				\item Transformation formulas (1):
				\begin{equation}\label{Pre 1.4}
				\begin{split}
				F(a,b;c;z)=(1-z)^{c-a-b} F(c-a,c-b;c;z).
				\end{split}
				\end{equation}
				
				\item Transformation formulas (2): if $c-a-b$ is not an integer, then
				\begin{equation}\label{Pre 1.5}
				\begin{split}
				F(a,b;c;z)=&\frac{\Gamma(c)\Gamma(c-a-b)}{\Gamma(c-a)\Gamma(c-b)} F(a,b;a+b-c+1;1-z)+\\
				&(1-z)^{c-a-b}\frac{\Gamma(c)\Gamma(a+b-c)}{\Gamma(a)\Gamma(b)}F(c-a,c-b;c-a-b+1;1-z).
				\end{split}
				\end{equation}
				
				\item Transformation formulas (3):
				\begin{align}\label{Pre 1.6}
				F\Big(a,b;2b;\frac{4z}{(1+z)^{2}}\Big)=(1+z)^{2a}F\Big(a,a+\frac{1}{2}-b;b+\frac{1}{2};z^{2}\Big).
				\end{align}
               \item if $\mathrm{Re} \mu>0$ and $|t|<1$, then (see \cite{Gradshteyn&Ryzhik}, page 407, 3.665)
			\begin{align}\label{Pre 1.7}
				\int_{0}^{\pi}\frac{\sin^{2\mu-1}\theta}{(1-2t\cos\theta+t^2)^{\alpha}}\ud \theta= \frac{\Gamma(\mu)\Gamma\left(\frac{1}{2}\right)}{\Gamma\left(\mu+\frac{1}{2}\right)}F\left(\alpha,\alpha-\mu+\frac{1}{2};\mu+\frac{1}{2}; t^2\right).
			\end{align}
			\end{itemize}
		
			\subsection{Funk-Hecke formula}

Let $K\in L^{1}((-1,1), (1-t^{2})^{(n-2)/2})$. The Funk-Hecke formula reads  (see \cite[Chapter 22]{Abramowitz&Stegun})
\begin{align}\label{Funk-Hecke}
			\int_{\S^n} K(\xi\cdot\eta) Y_l(\eta) \ud V_{{\S^n}}(\eta)=\lambda_lY_l(\xi),\quad Y_l \in \mathscr{H}_l,
			\end{align}
		where
			\begin{align*}
			\lambda_l=(4\pi)^{\frac{n-1}{2}}\frac{\Gamma\left(\frac{n-1}{2}\right)l!}{\Gamma\left(l+n-1\right)}\int_{-1}^{1}K(t)C^{\frac{n-1}{2}}_l(t)(1-t^2)^{\frac{n-2}{2}}\ud t
			\end{align*}
		and $ C_l^{\frac{n-1}{2}} $ denotes the $l$-th Gegenbauer polynomial with degree $\frac{n-1}{2}$.

In particular, we have that for $\gamma>0$, it holds (see \cite[(17)]{Beckner} or \cite[Corollary 4.3]{Frank&Lieb})
	\begin{align}\label{Funk-Hecke2}
	\int_{\S^n} |\xi-\eta|^{2\gamma-n} Y_l(\eta) \ud V_{{\S^n}}(\eta)=\frac{2^{2\gamma}\pi^{\frac{n}{2}}\Gamma(\gamma)}{\Gamma\left(\frac{n}{2}-\gamma\right)}
\frac{\Gamma\left(l+\frac{n}{2}-\gamma\right)}{\Gamma\left(l+\frac{n}{2}+\gamma\right)}Y_l(\xi),\quad Y_l \in \mathscr{H}_l.
	\end{align}

		\section{The fractional reverse Sobolev inequalities}\label{Sec 3}

	\subsection{First case: $\gamma\in \left(\frac{n}{2},\frac{n}{2}+1\right)$}

	We now present a notable sharp inequality involving the operators
	 $P_{2\gamma}$ and $P^{-1}_{2\gamma}$ , which serves as a central tool in establishing a reverse Sobolev inequality for values of $\gamma$ in the range
	  $\left(\frac{n}{2},\frac{n}{2}+1\right)$.
By (\ref{Funk-Hecke2}), the inverse operator $P^{-1}_{2\gamma}$ for $\gamma\notin \frac{n-2}{2}+\mathbb{N}$, can be defined  in terms of the integral kernel as follows:
	\begin{align}\label{3.1}
	P^{-1}_{2\gamma}(f)(\xi)=\frac{\Gamma\left(\frac{n}{2}-\gamma\right)}{2^{2\gamma}\pi^{\frac{n}{2}}\Gamma(\gamma)}\int_{\S^n} |\xi-\eta|^{2\gamma-n}f(\eta)  \ud V_{{\S^n}}(\eta).
	\end{align}
	The primary approach we use involves expanding positive functions on the sphere through spherical harmonics.
	
	\begin{lem}\label{Important Lem}
		Suppose $f\in H^{\gamma}(\S^n)$ and $g\in H^{-\gamma}(\S^n)$ are positive functions. If $\gamma\in \left(\frac{n}{2},\frac{n}{2}+1\right)$, then
		\begin{align}\label{Important formula 1}
		\int_{\S^n} f(-P_{2\gamma})(f)  \ud V_{{\S^n}}	\int_{\S^n} g(-P^{-1}_{2\gamma})(g)  \ud V_{{\S^n}}\leq \left(\int_{\S^n} fg  \ud V_{{\S^n}}\right)^2,
		\end{align}
		with the equality  holds if and only if $f=-cP_{2\gamma}^{-1}g$ with $c>0$.
	\end{lem}
	\begin{pf}
		Since $ \gamma \in \left( \frac{n}{2}, \frac{n}{2}+1 \right) $, it is easy to observe that $\Gamma\left( \frac{n}{2} - \gamma \right) < 0$. Therefore, we have
		\begin{align*}
		\int_{\S^n} g(-P^{-1}_{2\gamma})(g)  \ud V_{{\S^n}}
		=&-\frac{\Gamma\left(\frac{n}{2}-\gamma\right)}{2^{2\gamma}\pi^{\frac{n}{2}}\Gamma(\gamma)}\int_{\S^n} \int_{\S^n}g(\xi)g(\eta)|\xi-\eta|^{2\gamma-n}  \ud V_{{\S^n}}(\eta)\ud V_{{\S^n}}(\xi)>0.
		\end{align*}
		To prove \eqref{Important formula 1}, we  assume, without loss of generality, that
		\begin{align*}
		\int_{\S^n} f(-P_{2\gamma})(f)  \ud V_{{\S^n}}>0,
		\end{align*}
		since, in the case where this condition does not hold, the inequality in \eqref{Important formula 1} trivially holds. Now, we proceed by expanding the functions $f $ and $ g $ in terms of the spherical harmonics. Specifically, we write
		\begin{align*}
		f(\eta)=\sum_{l=0}^{+\infty}Y_l(f)(\eta) \qquad \mathrm{and}\qquad g(\eta)=\sum_{l=0}^{+\infty}Y_l(g)(\eta).
		\end{align*}
	For the sake of simplicity and clarity in notation, we introduce the following quantities
		\begin{align}\label{Important formula a}
		\int_{\S^n} f(-P_{2\gamma})(f)  \ud V_{{\S^n}}=&-\sum_{l=0}^{+\infty}\frac{\Gamma\left(l+\frac{n}{2}+\gamma\right)}{\Gamma\left(l+\frac{n}{2}-\gamma\right)}\int_{\S^n} |Y_l(f)|^2\ud V_{{\S^n}}\nonumber\\
		=:&a_f-|\mathbf{a}|^2>0,
		\end{align}
		where
		\begin{align*}
		a_f=-\frac{\Gamma\left(\frac{n}{2}+\gamma\right)}{\Gamma\left(\frac{n}{2}-\gamma\right)}|\S^n|Y^2_0(f) \qquad \mathrm{and}\qquad \mathbf{a}=\left(\sqrt{\frac{\Gamma\left(l+\frac{n}{2}+\gamma\right)}{\Gamma\left(l+\frac{n}{2}-\gamma\right)}\int_{\S^n} |Y_l(f)|^2\ud V_{{\S^n}}}\right)_{l\geq 1}.
		\end{align*}
		Similarly, for the function $ g $, we write
		\begin{align}\label{Important formula b}
		\int_{\S^n} g(-P^{-1}_{2\gamma})(g)  \ud V_{{\S^n}}=&-\sum_{l=0}^{+\infty}\frac{\Gamma\left(l+\frac{n}{2}-\gamma\right)}{\Gamma\left(l+\frac{n}{2}+\gamma\right)}\int_{\S^n} |Y_l(g)|^2\ud V_{{\S^n}}\nonumber\\
		=:&b_{g}-|\mathbf{b}|^2>0,
		\end{align}
		where
		\begin{align*}
		b_g=-\frac{\Gamma\left(\frac{n}{2}-\gamma\right)}{\Gamma\left(\frac{n}{2}+\gamma\right)}|\S^n|Y^2_0(g) \qquad \mathrm{and}\qquad \mathbf{b}=\left(\sqrt{\frac{\Gamma\left(l+\frac{n}{2}-\gamma\right)}{\Gamma\left(l+\frac{n}{2}+\gamma\right)}\int_{\S^n} |Y_l(g)|^2\ud V_{{\S^n}}}\right)_{l\geq 1}.
		\end{align*}
	Using the Cauchy-Schwarz inequality, it is straightforward to establish that
		\begin{align}\label{Important formula c}
		\int_{\S^n} fg  \ud V_{{\S^n}}=&Y_0(f)Y_0(g)|\S^n|+\sum_{l=1}^{+\infty}\int_{\S^n} Y_l(f)Y_{l}(g)  \ud V_{{\S^n}}\nonumber\\
		\geq &\sqrt{a_fb_{g}}-\mathbf{a}\cdot\mathbf{b}\geq \sqrt{a_fb_{g}}-|\mathbf{a}\|\mathbf{b}|>0.
		\end{align}
		Additionally, we observe a crucial equality
		\begin{align}\label{Important formula d}
		(a_f-|\mathbf{a}|^2)(b_{g}-|\mathbf{b}|^2)-(\sqrt{a_fb_{g}}-|\mathbf{a}\|\mathbf{b}|)^2=-(|\mathbf{a}|\sqrt{b_{g}}-|\mathbf{b}|\sqrt{a_f})^2\leq 0,
		\end{align}
	which leads to the final inequality following from \eqref{Important formula a}, \eqref{Important formula b}, \eqref{Important formula c} and \eqref{Important formula d}.

	When equality holds, we obtain the following identities for $ l \geq 1 $:
		\begin{align*}
		\int_{\S^n}Y_l(f)Y_l(g)   \ud V_{{\S^n}}=-\sqrt{\int_{\S^n} |Y_l(f)|^2\ud V_{{\S^n}}\int_{\S^n} |Y_l(g)|^2\ud V_{{\S^n}}} 	
		\end{align*}
		and
		\begin{align*}
		\mathbf{a}\cdot\mathbf{b}=|\mathbf{a}\|\mathbf{b}|,\qquad\qquad \frac{|\mathbf{a}|}{|\mathbf{b}|}=\frac{\sqrt{a_f}}{\sqrt{b_{g}}}=-
\frac{\Gamma\left(\frac{n}{2}+\gamma\right)}{\Gamma\left(\frac{n}{2}-\gamma\right)}\frac{\int_{\S^n} f  \ud V_{{\S^n}}}{\int_{\S^n}  g \ud V_{{\S^n}}}:=c_0.
		\end{align*}
	Thus, when $l\geq 1$,  we obtain $ |Y_l(f)| = c_l |Y_l(g)| $ for some nonnegative constant $c_l $,
		\begin{align*}
		Y_l(f)Y_l(g)\leq0 \qquad \mathrm{and}\qquad \frac{\sqrt{\int_{\S^n} |Y_l(f)|^2\ud V_{{\S^n}}}}{\sqrt{\int_{\S^n} |Y_l(g)|^2\ud V_{{\S^n}}}}\frac{\Gamma\left(l+\frac{n}{2}+\gamma\right)}{\Gamma\left(l+\frac{n}{2}-\gamma\right)}=c_0.
		\end{align*}
		So it  follows that
		\begin{align}\label{3.6}
		Y_l(f)=-c_0Y_l(g)\frac{\Gamma\left(l+\frac{n}{2}-\gamma\right)}{\Gamma\left(l+\frac{n}{2}+\gamma\right)} \qquad\mathrm{for}\qquad l\geq 1.
		\end{align}
Combing (\ref{3.6}) and $\sqrt{a_f}=c_0 \sqrt{b_{g}}$ yields
\begin{align*}
f=\sum_{l=0}^{+\infty}Y_l(f)=-c_0\frac{\Gamma\left(\frac{n}{2}-\gamma\right)}{\Gamma\left(\frac{n}{2}+\gamma\right)}Y_0(g)-c_0\sum_{l=1}^{+\infty}
\frac{\Gamma\left(l+\frac{n}{2}-\gamma\right)}{\Gamma\left(l+\frac{n}{2}+\gamma\right)}Y_l(g)=-c_0 P_{2\gamma}^{-1}g.
\end{align*}
The proof of Lemma \ref{Important Lem} is thereby completed.
	\end{pf}

Another crucial element in our proof is the Reverse Hardy-Littlewood-Sobolev inequality, an essential tool in analysis. This inequality, originally established by Dou and Zhu in \cite{Dou&Zhu}, provides a reverse form of the classical Hardy-Littlewood-Sobolev inequality, offering significant insights into the behavior of certain operators and integrals in functional spaces.
	\begin{thm}[Dou and Zhu]\label{Reverse HLS thm}
		For any nonnegative function $f,g\in L^{\frac{2n}{2n+\lambda}}(\S^n)$ with $\lambda>0$, it holds
		\begin{align*}
		\int_{\S^n}\int_{\S^n}&f(\xi)g(\eta)|\xi-\eta|^{\lambda}   \ud V_{{\S^n}}(\xi)   \ud V_{{\S^n}}(\eta)\\
		\geq&\pi^{-\frac{\lambda}{2}}\frac{\Gamma\left(\frac{n+\lambda}{2}\right)}{\Gamma\left(n+\frac{\lambda}{2}\right)}\left(\frac{\Gamma(n)}{\Gamma\left(\frac{n}{2}\right)}\right)^{1+\frac{\lambda}{n}}\|f\|_{L^{\frac{2n}{2n+\lambda}}(\S^n)}\|g\|_{L^{\frac{2n}{2n+\lambda}}(\S^n)}.
		\end{align*}
	Moreover, equality holds, up to a constant scaling, if and only if there exists a conformal transformation $\varphi$ such that  $\varphi^{*}g_{\S^n}=f^{\frac{4}{2n+\lambda}}g_{\S^n}$, i.e., $f=(\det\ud \varphi)^{\frac{2n+\lambda}{2n}}$.
	\end{thm}

	By applying Lemma \ref{Important Lem} and Theorem \ref{Reverse HLS thm}, we obtain the sharp reverse Sobolev inequality for $\gamma\in \left(\frac{n}{2},\frac{n}{2}+1\right)$.
	
	\begin{thm}\label{Sobolev inequ Case 1 thm}
		For any positive function $ f\in H^{\gamma}(\S^n)$ and $\gamma\in \left(\frac{n}{2},\frac{n}{2}+1\right)$, it holds
		\begin{align}\label{Sobolev inequ Case 1}
		\int_{\S^n}fP_{2\gamma}f   \ud V_{{\S^n}}\geq \frac{\Gamma\left(\frac{n}{2}+\gamma\right)}{\Gamma\left(\frac{n}{2}-\gamma\right)}|\S^n|^{\frac{2\gamma}{n}}\|f\|^2_{L^{\frac{2n}{n-2\gamma}}(\S^n)}.
		\end{align}
		Furthermore, equality is achieved, up to a constant scaling, if and only if there exists a conformal transformation
		 $\varphi$ such that $f=(\det\ud \varphi)^{\frac{n-2\gamma}{2n}}$.
	\end{thm}
	\begin{pf}
	First, by applying the reverse Hardy-Littlewood-Sobolev inequality and setting $\lambda = 2\gamma - n$, we obtain
		\begin{align*}
		\int_{\S^n} g(-P^{-1}_{2\gamma})(g)  \ud V_{{\S^n}}
		=&-\frac{\Gamma\left(\frac{n}{2}-\gamma\right)}{2^{2\gamma}\pi^{\frac{n}{2}}\Gamma(\gamma)}\int_{\S^n} \int_{\S^n}g(\xi)g(\eta)|\xi-\eta|^{2\gamma-n}  \ud V_{{\S^n}}(\eta)\ud V_{{\S^n}}(\xi)\\
		\geq&-\frac{\Gamma\left(\frac{n}{2}-\gamma\right)}{2^{2\gamma}\Gamma\left(\frac{n}{2}+\gamma\right)}\pi^{-\gamma}\left(\frac{\Gamma(n)}{\Gamma\left(\frac{n}{2}\right)}\right)^{\frac{2\gamma}{n}}\|g\|^2_{L^{\frac{2n}{n+2\gamma}}(\S^n)}.
		\end{align*}
		Substituting this inequality into \eqref{Important formula 1} and setting $g = f^{\frac{n + 2\gamma}{n - 2\gamma}}$, we derive that
		\begin{align*}
		-\frac{\Gamma\left(\frac{n}{2}-\gamma\right)}{2^{2\gamma}\Gamma\left(\frac{n}{2}+\gamma\right)}\pi^{-\gamma}\left(\frac{\Gamma(n)}{\Gamma\left(\frac{n}{2}\right)}\right)^{\frac{2\gamma}{n}}\int_{\S^n} f(-P_{2\gamma})(f)  \ud V_{{\S^n}}\leq \|f\|_{L^{\frac{2n}{n-2\gamma}}}^2.
		\end{align*}
		Recalling the Duplication formula
		\begin{align*}
		\Gamma\left(\frac{n}{2}\right)\Gamma\left(\frac{n+1}{2}\right)=2^{1-n}\sqrt{\pi}\Gamma(n)
		\end{align*}
		and using the identity $|\S^n| \Gamma\left(\frac{n+1}{2}\right) = 2 \pi^{\frac{n+1}{2}}$, we can obtain the desired Sobolev inequality \eqref{Sobolev inequ Case 1}. When equality holds, by Theorem \ref{Reverse HLS thm}, up to a constant, we have $f^{\frac{n + 2\gamma}{n - 2\gamma}} = g = (\det d\varphi)^{\frac{n + 2\gamma}{2n}}$. This completes  the proof of Theorem \ref{Sobolev inequ Case 1 thm}.
		
	\end{pf}

	\subsection{Second case: $\gamma\in \left(\frac{n}{2}+1,\frac{n}{2}+2\right)$}
	For $\gamma\in \left(\frac{n}{2}+1,\frac{n}{2}+2\right)$,we can provide a straightforward proof using a basic computation along with the reverse Holder inequality.
	\begin{thm}\label{Sobolev inequ Case 2 thm}
		For any positive function $ f\in H^{\gamma}(\S^n)$ and $\gamma\in \left(\frac{n}{2}+1,\frac{n}{2}+2\right)$, it holds
		\begin{align}\label{Sobolev inequ Case 2}
		\int_{\S^n}fP_{2\gamma}f   \ud V_{{\S^n}}\geq \frac{\Gamma\left(\frac{n}{2}+\gamma\right)}{\Gamma\left(\frac{n}{2}-\gamma\right)}|\S^n|^{\frac{2\gamma}{n}}\|f\|^2_{L^{\frac{2n}{n-2\gamma}}(\S^n)}.
		\end{align}
		Furthermore, equality holds, up to a constant scaling, if and only if there exists a conformal transformation
		 $\varphi$ such that $f=(\det\ud \varphi)^{\frac{n-2\gamma}{2n}}$.
	\end{thm}
	\begin{pf}
First, we aim to establish the inequality \eqref{Sobolev inequ Case 2} under the condition
		\begin{align*}
		\int_{\S^n} f(\xi)~\vec{\xi}\ud V_{{\S^n}}=0.
		\end{align*}
		To proceed, we express $ f $ in terms of its spherical harmonics expansion
		\begin{align*}
		f(\xi)=Y_0(f)+\sum_{l=2}^{+\infty}Y_l(f)(\xi) \qquad \mathrm{with}\qquad Y_0(f)=\fint_{\S^n}  f(\xi) \ud V_{{\S^n}}(\xi)>0.
		\end{align*}
		This decomposition allows us to rewrite the integral involving $ f $ as follows
		\begin{align}\label{Sobolev inequ Case 2 formula a0}
		\int_{\S^n}fP_{2\gamma}f   \ud V_{{\S^n}}=\frac{\Gamma\left(\frac{n}{2}+\gamma\right)}{\Gamma\left(\frac{n}{2}-\gamma\right)}|\S^n|Y^2_0(f)+\sum_{l=2}^{+\infty}\frac{\Gamma\left(l+\frac{n}{2}+\gamma\right)}{\Gamma\left(l+\frac{n}{2}-\gamma\right)}\int_{\S^n} |Y_l(f)|^2  \ud V_{{\S^n}}.
		\end{align}
		Since $\Gamma\left(\frac{n}{2} - \gamma\right) > 0$ and $\Gamma\left(l + \frac{n}{2} - \gamma\right) > 0$ for $l = 2, 3, \dots$, we can deduce that
		\begin{align}\label{Sobolev inequ Case 2 formula a}
		\int_{\S^n}fP_{2\gamma}f   \ud V_{{\S^n}}\geq \frac{\Gamma\left(\frac{n}{2}+\gamma\right)}{\Gamma\left(\frac{n}{2}-\gamma\right)}|\S^n|Y^2_0(f)=\frac{\Gamma\left(\frac{n}{2}+\gamma\right)}{\Gamma\left(\frac{n}{2}-\gamma\right)}|\S^n|^{-1}\left(\int_{\S^n} f  \ud V_{{\S^n}}\right)^2.
		\end{align}
	Next, applying the reverse Hölder inequality, we obtain
		\begin{align}\label{Sobolev inequ Case 2 formula b}
		\int_{\S^n}f\cdot 1   \ud V_{{\S^n}}\geq \|f\|_{L^{\frac{2n}{n-2\gamma}}}\left(\int_{\S^n}  \ud V_{{\S^n}}\right)^{\frac{n+2\gamma}{2n}}=\|f\|_{L^{\frac{2n}{n-2\gamma}}}|\S^n|^{\frac{n+2\gamma}{2n}}.
		\end{align}
	Substituting this inequality \eqref{Sobolev inequ Case 2 formula b} into \eqref{Sobolev inequ Case 2 formula a}, we get the desired inequality \eqref{Sobolev inequ Case 2}.

	Next, for a general function $ f \in H^{\gamma}(\S^n) $, it is a well-established fact that there exists a conformal transformation $ \varphi $ such that
(see Hang \cite{Hang})
		\begin{align*}
		\int_{\S^n}f_{{\varphi}}(\xi)~\vec{\xi}   \ud V_{{\S^n}}=0,
		\end{align*}
		where
		\begin{align*}
		f_{\varphi}(\xi)=f\circ\varphi(\xi)\left(\mathrm{\det}\ud \varphi(\xi)\right)^{\frac{n-2\gamma}{2n}}.
		\end{align*}
		Thanks to the conformal invariance of the operator $ P_{2\gamma} $ (see \cite[Corollary 2.2]{Hang} and \cite[Lemma 3]{Frank&Konig&Tang2}), we have the following equalities
		\begin{align*}
		\int_{\S^n} f_{\varphi}P_{2\gamma}(f_{\varphi})  \ud V_{{\S^n}}=\int_{\S^n}fP_{2\gamma}f   \ud V_{{\S^n}} \quad \mathrm{and}\quad \|f\|_{L^{\frac{2n}{n-2\gamma}}(\S^n)}=\|f_{\varphi}\|_{L^{\frac{2n}{n-2\gamma}}(\S^n)}.
		\end{align*}
Therefore,
\begin{align*}
\int_{\S^n}fP_{2\gamma}f   \ud V_{{\S^n}}=\int_{\S^n} f_{\varphi}P_{2\gamma}(f_{\varphi})  \ud V_{{\S^n}}\geq &
\frac{\Gamma\left(\frac{n}{2}+\gamma\right)}{\Gamma\left(\frac{n}{2}-\gamma\right)}|\S^n|^{\frac{2\gamma}{n}}\|f_{\varphi}\|^{2}_{L^{\frac{2n}{n-2\gamma}}(\S^n)}\\
=&\frac{\Gamma\left(\frac{n}{2}+\gamma\right)}{\Gamma\left(\frac{n}{2}-\gamma\right)}|\S^n|^{\frac{2\gamma}{n}}\|f\|^2_{L^{\frac{2n}{n-2\gamma}}(\S^n)}.
\end{align*}
		
		Finally, when equality holds in  (\ref{Sobolev inequ Case 2}), it follows from \eqref{Sobolev inequ Case 2 formula b} that $ f_{\varphi} \equiv \text{const} $. This implies that $ f $ must be of the form
		 $f=c\left(\mathrm{\det}\ud \varphi(\xi)\right)^{-\frac{n-2\gamma}{2n}}$ for some conformal transformation $\varphi$ on $\S^n$.
	\end{pf}

	\textbf{Proof of Theorem \ref{Sobolev Thm}:} By applying Theorem \ref{Sobolev inequ Case 1 thm} and Theorem \ref{Sobolev inequ Case 2 thm}, we have completed the proof.\\

Finally, we note that the reverse Sobolev inequalities also hold for nonnegative functions. In this regard, we present the following fundamental proposition concerning the Sobolev embedding theorem with a negative exponent, inspired by the work of Hang and Yang \cite{Hang&Yang1}.

		\begin{pro}\label{Positive Pro}
		Let $f$ be a nonnegative function with  $f\in H^{\gamma}(\S^n)$, where $\gamma\in \left(\frac{n}{2},\frac{n}{2}+1\right)\cup \left(\frac{n}{2}+1,\frac{n}{2}+2\right) $. If $\min_{\S^n}f=0$ then we have
		\begin{align*}
			\|f\|_{L^{\frac{2n}{n-2\gamma}}(\S^n)}=0 \qquad \mathrm{and}\qquad a_{2\gamma}(f)\geq 0.
		\end{align*}
	\end{pro}
	\begin{pf}
	For the first assertion, we prove this argument by two cases.
		
		\textbf{Case 1}: $
		\gamma\in \left(\frac{n}{2},\frac{n}{2}+1\right)$. By the Sobolev embedding theorem, we can get $H^{\gamma}(\S^n)\hookrightarrow C^{\gamma-\frac{n}{2}}(\S^n)$. Suppose that $f$ has a zero at some point on $\S^n$, without loss of generality, we can assume $f(S)=0$, where $S$ is the south pole. Letting $A=\|f\|_{C^{\gamma-\frac{n}{2}}(\S^n)}$, then
		\begin{align*}
		f(p)=|f(p)-f(S)|\leq Ar(p)^{\gamma-\frac{n}{2}},
		\end{align*}
		where $r(p)=\mathrm{dist}_{\S^n}(p,S)$. Noticing that $\frac{2n}{n-2\gamma}<0$,  we have
		\begin{align*}
		\int_{\S^n_{-}}f^{\frac{2n}{n-2\gamma}}\ud V_{\S^n}\geq A^{\frac{2n}{n-2\gamma}}\int_{\S^n_{-}}r^{-n}\ud V_{\S^n}=+\infty,
		\end{align*}
	where $\S^n_{-}=\mathbb{S}^{n}\cap\{(x_1,\cdots,x_{n+1})\in\mathbb{R}^{n+1}| x_{n+1}<0\}$.
		
		\textbf{Case 2}: $\gamma\in \left(\frac{n}{2}+1,\frac{n}{2}+2\right)$. In this case, we have  $H^{\gamma}(\S^n)\hookrightarrow C^{1,\gamma-\frac{n}{2}-1}(\S^n)$. Since $S$ is the minimizer of $f$,  we can get $\nabla_{\S^n}f(S)=0$. Letting $B=\|\nabla_{\S^n} f\|_{C^{\gamma-\frac{n}{2}-1}(\S^n)}$, one can  see
		\begin{align*}
		|\nabla_{\S^n} f(p)|=|\nabla_{\S^n} f(p)-\nabla_{\S^n} f (S)|\leq Br(p)^{\gamma-\frac{n}{2}-1},
		\end{align*}
		which leads to
		\begin{align*}
		f(p)=&|f(p)-f(S)|=\left|\int_{0}^{r(p)}\nabla_{\S^n} f(\sigma(t))\cdot\dot{\sigma}(t)\ud t\right|\\
		\leq& Br^{\gamma-\frac{n}{2}-1}(\sigma(t))r(p)\leq Br^{\gamma-\frac{n}{2}}(p)\qquad\qquad \qquad\mathrm{for}\qquad p\in \S^{n}_{-},
		\end{align*}
		where $\sigma(t)$ is the minimal geodesic from $S$ to $p$ on $\S^n$ with $|\dot{\sigma}(t)|=1$.
		The remaining part follows from the argument in the case $ \gamma\in \left(\frac{n}{2},\frac{n}{2}+1\right)$.
		
	For the second assertion, given any $\epsilon > 0$, by Theorem \ref{Sobolev Thm}, we have
		\begin{align}\label{Positive Pro formula a}
			a_{2\gamma}(f+\e)\geq \frac{\Gamma\left(\frac{n}{2}+\gamma\right)}{\Gamma\left(\frac{n}{2}-\gamma\right)}|\S^n|^{\frac{2\gamma}{n}}\|f+\e\|^2_{L^{\frac{2n}{n-2\gamma}}(\S^n)}.
		\end{align}
	By applying Fatou's Lemma, we obtain
		\begin{align*}
				\lim_{\e\to 0}\left(\int_{\S^n}(f+\e)^{\frac{2n}{n-2\gamma}}\ud V_{\S^n}\right)^{\frac{n-2\gamma}{2n}}=0.
		\end{align*}
		On the other hand, we also have
		\begin{align*}
				a_{2\gamma}(f+\e)=a_{2\gamma}(f)+\frac{\Gamma\left(\frac{n}{2}+\gamma\right)}{\Gamma\left(\frac{n}{2}-\gamma\right)}\left(2\e\int_{\S^n}f\ud V_{\S^n}+\e^2|\S^n|\right).
		\end{align*}
		Using this relation and taking the limit as \(\epsilon \to 0\) in \eqref{Positive Pro formula a}, we conclude that $a_{2\gamma}(f)\geq 0$.
		
	\end{pf}

		\textbf{Proof of Theorem \ref{Stability Thm}:}	
Let $ \varphi $ be the conformal transformation  such that
		\begin{align*}
		\int_{\S^n}f_{{\varphi}}(\xi)~\vec{\xi}   \ud V_{{\S^n}}=0,
		\end{align*}
where $f_{\varphi}(\xi)=f\circ\varphi(\xi)\left(\mathrm{\det}\ud \varphi(\xi)\right)^{\frac{n-2\gamma}{2n}}.$
We note that $f_{\varphi}$ has the expansion in spherical harmonics
$$f_{\varphi}=Y_0(f_{\varphi})+\sum\limits_{l=2}^{+\infty}Y_l(f_{\varphi})\quad \mathrm{with}\quad Y_0(f_{\varphi})=\fint_{\S^n}  f_{\varphi}(\xi) \ud V_{{\S^n}}(\xi).$$
We compute, by using (\ref{Sobolev inequ Case 2 formula b}),
\begin{align*}
		a_{2\gamma}(f_{{\varphi}})=&\frac{\Gamma\left(\frac{n}{2}+\gamma\right)}{\Gamma\left(\frac{n}{2}-\gamma\right)}|\S^n|Y^2_0(f_{{\varphi}})+
\sum_{l=2}^{+\infty}\frac{\Gamma\left(l+\frac{n}{2}+\gamma\right)}{\Gamma\left(l+\frac{n}{2}-\gamma\right)}\int_{\S^n} |Y_l(f_{{\varphi}})|^2  \ud V_{{\S^n}}\\
		=&\frac{\Gamma\left(\frac{n}{2}+\gamma\right)}{\Gamma\left(\frac{n}{2}-\gamma\right)}|\S^n|Y^2_0(f_{{\varphi}})+a_{2\gamma}(f_{{\varphi}}-\bar{f_{{\varphi}}})\\
		\geq&\frac{\Gamma\left(\frac{n}{2}+\gamma\right)}{\Gamma\left(\frac{n}{2}-\gamma\right)}|\S^n|^{\frac{2\gamma}{n}}\|f_{\varphi}\|^2_{L^{\frac{2n}{n-2\gamma}}(\S^n)}+
a_{2\gamma}(f_{{\varphi}}-\bar{f_{{\varphi}}}),
		\end{align*}
where $\bar{f_{\varphi}}=Y_0(f_{\varphi})=\fint_{\S^n} f_{\varphi}  \ud V_{{\S^n}}$. Obviously,
\begin{align*}
 a_{2\gamma}(f_{{\varphi}}-\bar{f_{{\varphi}}})=\sum_{l=2}^{+\infty}\frac{\Gamma\left(l+\frac{n}{2}+\gamma\right)}{\Gamma\left(l+\frac{n}{2}-\gamma\right)}\int_{\S^n} |Y_l(f_{{\varphi}})|^2  \ud V_{{\S^n}}\geq 0,
\end{align*}
with equality if and only if $f_{\varphi}=\bar{f_{\varphi}}$, i.e., $f$ is a extremal function of reverse Sobolev inequality.
Therefore, if $f$ is not a extremal function of reverse Sobolev inequality, then we have
\begin{align*}
a_{2\gamma}(f_{{\varphi}}-\bar{f_{{\varphi}}})>0
\end{align*}
and
\begin{align}\nonumber
 a_{2\gamma}(f)=a_{2\gamma}(f_{{\varphi}})
 \geq&\frac{\Gamma\left(\frac{n}{2}+\gamma\right)}{\Gamma\left(\frac{n}{2}-\gamma\right)}|\S^n|^{\frac{2\gamma}{n}}\|f_{\varphi}\|^2_{L^{\frac{2n}{n-2\gamma}}(\S^n)}+
a_{2\gamma}(f_{{\varphi}}-\bar{f_{{\varphi}}})\\
\label{3.13}
=&\frac{\Gamma\left(\frac{n}{2}+\gamma\right)}{\Gamma\left(\frac{n}{2}-\gamma\right)}|\S^n|^{\frac{2\gamma}{n}}\|f\|^2_{L^{\frac{2n}{n-2\gamma}}(\S^n)}+
a_{2\gamma}(f_{{\varphi}}-\bar{f_{{\varphi}}}).
\end{align}
	By the conformal invariance of $a_{2\gamma}$, we have
		\begin{align}\label{3.14}
		a_{2\gamma}(f_{\varphi}-\bar{f_{\varphi}})= a_{2\gamma}\left(f-\bar{f_{\varphi}}\left(\mathrm{\det}\ud \varphi^{-1}\right)^{\frac{n-2\gamma}{2n}}\right)>0.
		\end{align}
Substituting (\ref{3.14}) into (\ref{3.13}), we obtain
	\begin{align*}
   a_{2\gamma}(f)\geq \frac{\Gamma\left(\frac{n}{2}+\gamma\right)}{\Gamma\left(\frac{n}{2}-\gamma\right)}|\S^n|^{\frac{2\gamma}{n}}\|f\|^2_{L^{\frac{2n}{n-2\gamma}}(\S^n)}+
a_{2\gamma}\left(f-\bar{f_{\varphi}}\left(\mathrm{\det}\ud \varphi^{-1}\right)^{\frac{n-2\gamma}{2n}}\right).
 \end{align*}
This completes the proof of Theorem \ref{Stability Thm}.

	\section{Sobolev trace inequality}\label{Sec 4}
	\subsection{Scattering theory in the Poincar\'e ball}
In  this subsection, we will present the explicit integral formula for the solution of scattering theory within the context of the Poincar\'e  ball.
We  refer to \cite{ma,Graham&Zworski,Chang&Gonzalez} and the references therein for the scattering theory on manifolds.
Utilizing this integral formula, we will subsequently derive the series expansion with respect to
 $\rho$. This series expansion is of significant importance as it will facilitate our further analysis and enable us to derive the boundary operators effectively.
	\begin{thm}\label{Scattering thm}
	Let $\gamma \in (0, +\infty) \setminus \mathbb{N}$ and $s = \frac{n}{2} + \gamma$. For $f \in C(\S^n)$, there exists a unique function $u\in C^2(\B^{n+1})$ that solves the Poisson problem
		\begin{align}\label{Sec 5 equ}
			\begin{cases}
				\displaystyle -\Delta_{+}u-s(n-s)u=0 \qquad&\mathrm{in}\qquad \B^{n+1},\\
				\displaystyle \lim_{\rho_0\to 0^{+}}\rho_0^{s-n}u=f \qquad&\mathrm{on}\qquad \S^n.
			\end{cases}
		\end{align}
		The integral formula for the solution is given by
		\begin{align}\label{Sec 5 Integral formu solution}
			u(x)=\pi^{-\frac{n}{2}}\frac{2^{-s}\Gamma\left(\frac{n}{2}+\gamma\right)}{\Gamma(\gamma)}\int_{\S^n}  \left(\frac{1-|x|^2}{|x-\xi|^2}\right)^sf(\xi) \ud V_{{\S^n}}(\xi).
		\end{align}
		Furthermore, if $f$ has a spherical harmonic expansion of the form $f = \sum\limits_{l=0}^{+\infty} Y_l$ with $Y_l\in \mathscr{H}_l$, then the solution $u(x)$ takes the series form
	\begin{align}\label{Sec 5 Series formu solution}
		u(x)=\rho_0^{n-s}\sum_{l=0}^{+\infty}\varphi_l(r^2)r^{l}Y_l,
	\end{align}
where
	\begin{align}\label{Sec 5 Series varphi}
		\varphi_l(r)=&\frac{F\left(l+\frac{n}{2}-\gamma,\frac{1}{2}-\gamma; l+\frac{n+1}{2}; r\right)}{F\left(l+\frac{n}{2}-\gamma,\frac{1}{2}-\gamma; l+\frac{n+1}{2}; 1\right)}\nonumber\\
		=&\frac{\Gamma\left(\gamma+\frac{1}{2}\right)}{\Gamma(2\gamma)}\frac{\Gamma\left(l+\gamma+\frac{n}{2}\right)}{\Gamma(l+\frac{n+1}{2})}F\left(l+\frac{n}{2}-\gamma,\frac{1}{2}-\gamma; l+\frac{n+1}{2}; r\right).
	\end{align}
	\end{thm}

	\begin{pf}
		
		\textbf{Step 1}:
	We begin by verifying that the given integral formula satisfies the interior equation in $ \mathbb{B}^{n+1}$, specifically,
		\begin{align*}
			-\Delta_{+}u-s(n-s)u=0.
		\end{align*}
It is enough to show
		\begin{align}\label{Sec 5 thm 1 formula a}
			-\Delta_{+}\left(\frac{1-|x|^2}{|x-\xi|^2}\right)^{s}-s(n-s)\left(\frac{1-|x|^2}{|x-\xi|^2}\right)^{s}=0.
		\end{align}
		To prove the  identity \eqref{Sec 5 thm 1 formula a}, we first introduce the conformal Laplacian in $(\B^{n+1},g_{\B})$,
		\begin{align*}
			L_{g_{\B}}(u)=-\Delta_{+}u-\frac{n^2-1}{4}u.
		\end{align*}
	Through a straightforward calculation, 	then we can demonstrate that the above equation \eqref{Sec 5 thm 1 formula a} transforms into
		\begin{align}\label{Sec 5 thm 1 formula b}
			L_{g_{\B}}\left(\frac{1-|x|^2}{|x-\xi|^2}\right)^{s}=\left(\frac{1-(n-2s)^2}{4}\right)\left(\frac{1-|x|^2}{|x-\xi|^2}\right)^{s}.
		\end{align}
		Next, we rewrite the metric $g_{\B}=\rho_0^{-2}|\ud x|^2$ in the form
		$$g_{\B}=\left(\rho_0^{-\frac{n-1}{2}}\right)^{\frac{4}{n-1}}|\ud x|^2$$
	and, using the conformal invariance, we conclude that the above equation \eqref{Sec 5 thm 1 formula b} is equivalent to	
		\begin{align}\label{Sec 5 thm 1 formula c}
			-\Delta_{\R^{n+1}}\left[\rho_0^{-\frac{n-1}{2}}\left(\frac{1-|x|^2}{|x-\xi|^2}\right)^{s}\right]=\left(\frac{1-(n-2s)^2}{4}\right)\rho_0^{-\frac{n+3}{2}}\left(\frac{1-|x|^2}{|x-\xi|^2}\right)^{s}.
		\end{align}
		The final equation \eqref{Sec 5 thm 1 formula c} can be reduced to the following three identities:
		\begin{align*}
			\Delta_{\R^{n+1}}(1-|x|^2)^{s-\frac{n-1}{2}}=&-2(n+1)\left(s-\frac{n-1}{2}\right)(1-|x|^2)^{s-\frac{n+1}{2}}\\
			&+4\left(s-\frac{n-1}{2}\right)\left(s-\frac{n+1}{2}\right)|x|^2(1-|x|^2)^{s-\frac{n+3}{2}},\\
			\nabla (1-|x|^2)^{s-\frac{n-1}{2}}\cdot \nabla|x-\xi|^{-2s}=&4s\left(s-\frac{n-1}{2}\right)(1-|x|^2)^{s-\frac{n+1}{2}}|x-\xi|^{-2-2s}\langle x,x-\xi\rangle
		\end{align*}
		and
		\begin{align*}
			\Delta_{\R^{n+1}}|x-\xi|^{-2s}=-2s(n-2s-1)|x-\xi|^{-2s-2}.
		\end{align*}
		
		\textbf{Step 2:}
	 Verifying the boundary condition, namely, we need to prove the following limitation
		\begin{align}\label{Sec 5 thm 1 boundary data}
		\lim_{\rho_0\to 0^{+}}\rho_0^{s-n}u=f.
		\end{align}
		Consider $x_0\in \B^{n+1}$ and set $\xi_0=x_0/|x_0|\in \S^n$. For any small $\e>0$, there exists $\delta=\delta(\e)$ such that
		\begin{align*}
				|f(\xi)-f(\xi_0)|<\e\qquad\mathrm{for}\qquad \xi\in B_{\delta}(\xi_0),
		\end{align*}
	where we assume $|x_0| $ is near the boundary $\S^n$ and $1-|x_0|<\delta$. We can estimate the integral as follows:
		\begin{align}\label{Sec 5 thm 1 formula d}
			&\int_{\S^n} \frac{f(\xi)}{|x_0-\xi|^{2s}}  \ud V_{{\S^n}}(\xi)
			-f(\xi_0)\int_{\S^n}\frac{1}{|x_0-\xi|^{2s}}  \ud V_{{\S^n}}(\xi)\nonumber\\
			=&\int_{B_{\delta}(\xi_0)}\frac{f(\xi)-f(\xi_0)}{|x_0-\xi|^{2s}}  \ud V_{{\S^n}}(\xi)
			+\int_{\S^n\backslash B_{\delta}(\xi_0)}\frac{f(\xi)-f(\xi_0)}{|x_0-\xi|^{2s}}  \ud V_{{\S^n}}(\xi)\nonumber\\
			=&O(\e)	\int_{\S^n} \frac{1}{|x_0-\xi|^{2s}}  \ud V_{{\S^n}}(\xi)+O(\delta^{-2s}).
		\end{align}
		Using the Funk-Hecke formula \eqref{Funk-Hecke}, we have
		\begin{align*}
			\int_{\S^n}\frac{1}{|x_0-\xi|^{2s}}  \ud V_{{\S^n}}(\xi)=&\int_{\S^n}  \frac{1}{(|x_0|^2+1-2|x_0|\xi_0\cdot \xi)^{s}}  \ud V_{{\S^n}}(\xi)\\
			=&|\S^{n-1}|\int_{-1}^{1}\frac{1}{(|x_0|^2+1-2|x_0|t)^{s}} (1-t^2)^{\frac{n-2}{2}}\ud t\;\;\;(\textrm{substituting }\; t=\cos\theta)\\
			=&|\S^{n-1}|\int_{0}^{\pi}\frac{\sin^{n-1}\theta}{(|x_0|^2+1-2|x_0|\cos\theta)^{s}} \ud \theta.
		\end{align*}
		Then, applying \eqref{Pre 1.7},  \eqref{Pre 1.4} and \eqref{Pre 1.3}, we conclude that
		\begin{align}\label{Sec 5 thm 1 formula e}
			\int_{\S^n} \frac{1}{|x_0-\xi|^{2s}}  \ud V_{{\S^n}}(\xi)=&|\S^{n-1}|\frac{\Gamma\left(\frac{n}{2}\right)\Gamma\left(\frac{1}{2}\right)}{\Gamma\left(\frac{n+1}{2}\right)}F\left(s,s-\frac{n}{2}+\frac{1}{2}; \frac{n+1}{2}; |x_0|^2\right)\nonumber\\
			=&|\S^n|(1-|x_0|^2)^{n-2s}F\left(\frac{n+1}{2}-s,n-s; \frac{n+1}{2}; |x_0|^2\right)\nonumber\\
			=&|\S^n|(1-|x_0|^2)^{n-2s}\frac{\Gamma\left(\frac{n+1}{2}\right)\Gamma(2s-n)}{\Gamma(s)\Gamma\left(s-\frac{n-1}{2}\right)}(1+o(1)),
		\end{align}
		where $o(1)$ means that $o(1)\to 0$ as $|x_0|\to 1$. Using the formula \eqref{Sec 5 thm 1 formula d}, \eqref{Sec 5 thm 1 formula e}  and Duplication formula
 for Gamma function
	\begin{align}\label{Duplic formula}
	\Gamma(z)\Gamma\left(z+\frac{1}{2}\right)=2^{1-2z}\sqrt{\pi}\Gamma(2z),
	\end{align}
 we get
		\begin{align*}
		&	\left|\rho_0^{s-n}\pi^{-\frac{n}{2}}\frac{2^{-s}\Gamma\left(\frac{n}{2}+\gamma\right)}{\Gamma(\gamma)}\int_{\S^n}  \left(\frac{1-|x_0|^2}{|x_0-\xi|^2}\right)^s f(\xi)\ud V_{{\S^n}}(\xi)-f(\xi_0)\right|\\
			&\qquad\leq o(1)+C\e+C(1-|x_0|)^{2s-n}\delta^{-2s}.
		\end{align*}
	Hence, we conclude that for any $\e>0$, there holds
		\begin{align*}
			\limsup_{|x|\to 1}\left|\rho_0^{s-n}u(x)-f(x/|x|)\right|\leq C\e,
		\end{align*}
		which completes the proof of \eqref{Sec 5 thm 1 boundary data}.
		
		\textbf{Step 3:} Derivation of formula \eqref{Sec 5 Series formu solution} and uniqueness.
		Choosing an orthogonal basis $Y_{l,1}, Y_{l,2},\cdots, Y_{l, N_l}$ of $\mathscr{H}_{l}$ with $\|Y_{l, i}\|_{L^{2}(\mathbb{S}^{n})}=1(1\leq i\leq N_l)$,
we can write
$u(x)$ as spherical decomposition
		\begin{align*}
			u(x)=\sum_{l=0}^{+\infty}\sum_{i=1}^{N_l}u_{l,i}(r)Y_{l,i}(\xi),\;\; 0<r=|x|<1,
		\end{align*}
		where
		\begin{align}\label{Sec 5 thm 1 formula e0}
			u_{l,i}(r)=\int_{\S^n}  u(r\xi)Y_{l,i}(\xi) \ud V_{{\S^n}}(\xi) \quad \mathrm{and}\quad \lim_{r\to 0^{+}}u_{l,i}(r)=\begin{cases}
			\displaystyle |\mathbb{S}^{n}|^{1/2}u(0) \quad&\mathrm{for}\quad l=0,\\
			\displaystyle 0 \quad&\mathrm{for}\quad l\geq 1.
			\end{cases}
		\end{align}
		Using the Step 1, we find
		\begin{align*}
			-\Delta_{+}(u_{l,i}(r)Y_{l,i})-s(n-s)u_{l,i}(r)Y_{l,i}=0.
		\end{align*}
		Setting $v=\rho_0^{s-n}u$ and
		\begin{align*}
				v(x)=\sum_{l=0}^{+\infty}\sum_{i=1}^{N_l}v_{l,i}(r^{2})r^{l}Y_{l,i}(\xi),
		\end{align*}
		it is not hard to see that
		\begin{align}\label{Sec 5 thm 1 formula e00}
			u_{l,i}(r)=v_{l,i}(r^2)\rho_0^{n-s}r^{l}.
		\end{align}
		Following the similar techniques in \eqref{Sec 5 thm 1 formula c}, we have
		\begin{align*}
			-\Delta_{\R^{n+1}}\left(\rho_0^{\frac{n+1}{2}-s}v_{l,i}(r^2)r^{l}Y_{l,i}\right)-\left(\frac{1}{4}-\gamma^2\right)\rho_0^{\frac{n-3}{2}-s}v_{l,i}(r^2)r^{l}Y_{l,i}=0.
		\end{align*}
		Using the formulas $-\Delta_{\S^n}Y_{l,i}=l(n+l-1)Y_{l,i}$ and
		\begin{align*}
		\Delta_{\R^{n+1}}=\frac{\partial^2}{\partial r^2}+\frac{n}{r}\frac{\partial}{\partial r}+\frac{1}{r^2}\Delta_{\S^n},
		\end{align*}
		we have
		\begin{align}\label{Sec 5 thm 1 formula f}
			-\left(\frac{\partial^2}{\partial r^2}+\frac{n}{r}\frac{\partial}{\partial r}\right)\left(\rho_0^{\frac{n+1}{2}-s}v_{l,i}(r^2)r^{l}\right)+&l(n+l-1)\rho_0^{\frac{n+1}{2}-s}v_{l,i}(r^2)r^{l-2}\nonumber\\
			=&\left(\frac{1}{4}-\gamma^2\right)\rho_0^{\frac{n-3}{2}-s}v_{l,i}(r^2)r^{l}.
		\end{align}
		By changing the variable $t=r^2$, we obtain
		\begin{align}\label{Sec 5 thm 1 formula g}
			\rho_0=\frac{1-t}{2},\qquad \frac{\partial}{\partial r}=2r\frac{\partial}{\partial t}.
		\end{align}
		So it is not hard to see that
		\begin{align}\label{Sec 5 thm 1 formula i}
			\frac{\partial}{\partial r}\left(\rho_0^{\frac{n+1}{2}-s}v_{l,i}(r^2)r^{l}\right)
			=&\rho_0^{\frac{n-1}{2}-s}r^{l-1}\left[\left(l\rho_0+\left(\gamma-\frac{1}{2}\right)r^2\right)v_{l,i}(t)+2\rho_0r^2v^{'}_{l,i}(t)\right]
		\end{align}
		and
		\begin{align}\label{Sec 5 thm 1 formula j}
			&\frac{\partial^2}{\partial r^2}\left(\rho_0^{\frac{n+1}{2}-s}v_{l,i}(r^2)r^{l}\right)\nonumber\\
		=&\rho_0^{\frac{n-3}{2}-s}r^{l-2}\left[\left(\gamma+\frac{1}{2}\right)r^2+(l-1)\rho_0\right]\left[\left(l\rho_0+\left(\gamma-\frac{1}{2}\right)r^2\right)v_{l,i}(t)+2\rho_0r^2v^{'}_{l,i}(t)\right]\nonumber\\
		+&\rho_0^{\frac{n-1}{2}-s}r^{l}\left[(2\gamma-l-1)v_{l,i}(t)+(2(l+2)\rho_0+(2\gamma-3)r^2)v^{'}_{l,i}(t)+4\rho_0r^2v^{''}_{l,i}(t)\right].
		\end{align}
	By using identities \eqref{Sec 5 thm 1 formula f}, \eqref{Sec 5 thm 1 formula g}, \eqref{Sec 5 thm 1 formula i} and \eqref{Sec 5 thm 1 formula j}, we conclude that $v_{l,i}$ satisfies the hypergeometric equation
		\begin{align}\label{Sec 5 thm 1 formula k}
			t(1-t)v^{''}_{l,i}(t)+&\left[l+\frac{n+1}{2}-\left(l+\frac{n+3}{2}-2\gamma\right)t\right]v^{'}_{l,i}(t)
			\nonumber\\
			-&\left(l+\frac{n}{2}-\gamma\right)\left(\frac{1}{2}-\gamma\right)v_{l,i}(t)=0.
		\end{align}
		From \eqref{Sec 5 thm 1 formula e0} and \eqref{Sec 5 thm 1 formula e00}, we have the initial data
		\begin{align*}
			\lim_{t\to 0^{+}}t^{l/2}|v_{l,i}(t)| <+\infty.
		\end{align*}
		According to the properties of the hypergeometric equation (see page 1011 in \cite{Gradshteyn&Ryzhik}),  there exist two linear independent solutions $\omega_1(t)$, $\omega_2(t)$, where
		\begin{align*}
		\omega_1(t)=F\left(l+\frac{n}{2}-\gamma, \frac{1}{2}-\gamma; l+\frac{n+1}{2}; t\right)
		\end{align*}
		and $\omega_2(t)$ is a singular solution with asymptotic behaviour
		\begin{align*}
		\omega_2(t)\sim\begin{cases}
		\displaystyle \ln t  &n=1, l=0\\
		\displaystyle t^{{1-\frac{n+1}{2}-l}}\qquad&\mathrm{other}
		\end{cases}
	\qquad	\mathrm{as} \qquad t\to 0_{+}.
		\end{align*}
	In either case, we can get
	\begin{align*}
		\lim_{t\to 0^{+}}t^{l/2}|\omega_2(t)|=\infty,
	\end{align*}
which implies  $v_{l,i}(t)=c_{l,i}\omega_1(t)$.
Therefore, we get
\begin{align*}
			u(x)=\rho_0^{n-s}\sum_{l=0}^{+\infty}F\left(l+\frac{n}{2}-\gamma, \frac{1}{2}-\gamma; l+\frac{n+1}{2}; r^{2}\right)
r^{l}\left(\sum_{i=1}^{N_l}c_{l,i}Y_{l,i}(\xi)\right),\;\; 0<r=|x|<1.
	\end{align*}
Furthermore, since $	\lim\limits_{\rho_0\to 0^{+}}\rho_0^{s-n}u=f=\sum\limits_{l=0}^{+\infty} Y_l$, we find
	\begin{align*}
	Y_{l}=&\lim_{r\to 1^{-}}F\left(l+\frac{n}{2}-\gamma, \frac{1}{2}-\gamma; l+\frac{n+1}{2}; r^{2}\right)r^{l}
\left(\sum_{i=1}^{N_l}c_{l,i}Y_{l,i}(\xi)\right)\\
	=&\frac{\Gamma\left(l+\frac{n+1}{2}\right)\Gamma(2\gamma)}{\Gamma\left(\gamma+\frac{1}{2}\right)\Gamma\left(l+\gamma+\frac{n}{2}\right)}
\left(\sum_{i=1}^{N_l}c_{l,i}Y_{l,i}(\xi)\right).
	\end{align*}	i.e.,
\begin{align*}
\sum_{i=1}^{N_l}c_{l,i}Y_{l,i}(\xi)=\frac{\Gamma\left(\gamma+\frac{1}{2}\right)\Gamma\left(l+\gamma+\frac{n}{2}\right)}
{\Gamma\left(l+\frac{n+1}{2}\right)\Gamma(2\gamma)}
Y_l.
\end{align*}
Thus, we obtain the formula \eqref{Sec 5 Series formu solution}. The uniqueness part follows by taking $f\equiv 0$, which implies $c_{l,i}=0$.
	\end{pf}
	
	In the theorem above, we can express
	 $u$ as a series in terms of
	  $\rho_0$. However, we aim to derive an explicit formula in terms of
	   $\rho$, which requires the following calculations.
	\begin{lem}
		For $\gamma\in (0,+\infty)\backslash \mathbb{N}$ and $s=\frac{n}{2}+\gamma$, the solution of \eqref{Sec 5 equ} can be written as
		\begin{align}\label{Sec 5 lem 1 formula}
			u(x)=&\rho^{n-s}\sum_{l=0}^{+\infty}\left(1-\frac{\rho^2}{4}\right)^{l}F\left(l+\frac{n}{2}-\gamma, l+\frac{n}{2}; 1-\gamma; \frac{\rho^2}{4}\right)Y_l+\nonumber\\
			&\rho^{s}\frac{\Gamma(-\gamma)}{2^{2\gamma}\Gamma(\gamma)}\sum_{l=0}^{+\infty}\frac{\Gamma\left(l+\frac{n}{2}+\gamma\right)}{\Gamma\left(l+\frac{n}{2}-\gamma\right)}
\left(1-\frac{\rho^2}{4}\right)^{l}F\left(l+\frac{n}{2}+\gamma, l+\frac{n}{2}; 1+\gamma; \frac{\rho^2}{4}\right)Y_l.
		\end{align}
	\end{lem}
	\begin{pf}
First, we assume that $\gamma\in (0,+\infty)\backslash\frac{1}{2}\mathbb{N}$.		From equations \eqref{Pre 1.4} and \eqref{Pre 1.5}, we obtain the relationship
		\begin{align*}
			&F\left(l+\frac{n}{2}-\gamma,\frac{1}{2}-\gamma; l+\frac{n+1}{2}; r^2\right)\\
=&(1-r^2)^{2\gamma}F\left(\frac{1}{2}+\gamma,l+\frac{n}{2}+\gamma; l+\frac{n+1}{2}; r^2\right)\\
		=&(1-r^2)^{2\gamma}\frac{\Gamma\left(l+\frac{n+1}{2}\right)\Gamma(-2\gamma)}{\Gamma\left(l+\frac{n}{2}-\gamma\right)\Gamma\left(\frac{1}{2}-\gamma\right)}	F\left(\frac{1}{2}+\gamma,l+\frac{n}{2}+\gamma; 1+2\gamma; 1-r^2\right)+\\		&\frac{\Gamma\left(l+\frac{n+1}{2}\right)\Gamma(2\gamma)}{\Gamma\left(l+\frac{n}{2}+\gamma\right)\Gamma\left(\frac{1}{2}+\gamma\right)}F\left(l+\frac{n}{2}-\gamma, \frac{1}{2}-\gamma; 1-2\gamma; 1-r^2\right).
		\end{align*}
	Thus, from the formula \eqref{Sec 5 Series varphi}, we can express $\varphi_l(r^2)$ as follows:
		\begin{align*}
\varphi_l(r^2)=			&\left(\frac{2\rho}{\left(1+\frac{\rho}{2}\right)^2}\right)^{2\gamma}\frac{\Gamma\left(\gamma+\frac{1}{2}\right)}{\Gamma(2\gamma)}\frac{\Gamma\left(l+\frac{n}{2}+\gamma\right)\Gamma(-2\gamma)}{\Gamma(l+\frac{n}{2}-\gamma)\Gamma(\frac{1}{2}-\gamma)}
			F\left(\frac{1}{2}+\gamma, l+\frac{n}{2}+\gamma; 1+2\gamma; \frac{2\rho}{\left(1+\frac{\rho}{2}\right)^2}\right)\\
			&\qquad\qquad\qquad+F\left(l+\frac{n}{2}-\gamma, \frac{1}{2}-\gamma; 1-2\gamma; \frac{2\rho}{\left(1+\frac{\rho}{2}\right)^2}\right).
		\end{align*}
		Using the  Duplication formula
		\begin{align*}
			\Gamma(2\gamma)=2^{2\gamma-1}\frac{\Gamma(\gamma)\Gamma\left(\gamma+\frac{1}{2}\right)}{\Gamma\left(\frac{1}{2}\right)}, \qquad 	\Gamma(-2\gamma)=2^{-2\gamma-1}\frac{\Gamma(-\gamma)\Gamma\left(-\gamma+\frac{1}{2}\right)}{\Gamma\left(\frac{1}{2}\right)},
		\end{align*}
		and combining this with \eqref{Pre 1.6}, we get
		\begin{align}\label{Sec 5 lem 1 formula a}
			\varphi_l(r^2)=&\left(\frac{\rho}{\left(1+\frac{\rho}{2}\right)^2}\right)^{2\gamma}\frac{\Gamma\left(-\gamma\right)}{2^{2\gamma}\Gamma(\gamma)}\frac{\Gamma\left(l+\frac{n}{2}+\gamma\right)}{\Gamma(l+\frac{n}{2}-\gamma)}
			F\left(\frac{1}{2}+\gamma, l+\frac{n}{2}+\gamma; 1+2\gamma; \frac{2\rho}{\left(1+\frac{\rho}{2}\right)^2}\right)\nonumber\\
			&\qquad\qquad\qquad+F\left(l+\frac{n}{2}-\gamma, \frac{1}{2}-\gamma; 1-2\gamma; \frac{2\rho}{\left(1+\frac{\rho}{2}\right)^2}\right)\nonumber\\
			=&\rho^{2\gamma}\left(1+\frac{\rho}{2}\right)^{2l+n-2\gamma}\frac{\Gamma\left(-\gamma\right)}{2^{2\gamma}\Gamma(\gamma)}\frac{\Gamma\left(l+\frac{n}{2}+\gamma\right)}{\Gamma(l+\frac{n}{2}-\gamma)}
			F\left(l+\frac{n}{2}+\gamma, l+\frac{n}{2}; 1+\gamma; \frac{\rho^2}{4}\right)\nonumber\\
			&\qquad\qquad\qquad+\left(1+\frac{\rho}{2}\right)^{2l+n-2\gamma}F\left(l+\frac{n}{2}-\gamma, l+\frac{n}{2}; 1-\gamma; \frac{\rho^2}{4}\right).
		\end{align}
		By substituting $r=\frac{2-\rho}{2+\rho}$ and $\rho_0=\frac{\rho}{\left(1+\rho/2\right)^2}$ into equations \eqref{Sec 5 Series formu solution} and \eqref{Sec 5 lem 1 formula a}, we derive the  formula \eqref{Sec 5 lem 1 formula} for $\gamma\in (0,+\infty)\backslash\frac{1}{2}\mathbb{N}$.

When $\gamma\in \frac{1}{2}\mathbb{N}\backslash\mathbb{N}$, taking limit for the above formula \eqref{Sec 5 lem 1 formula} as $\gamma+\e\to \gamma$.
	\end{pf}
	\begin{rem}\label{Vital Rem}
		To be consistent with the scattering theory, we  let $s=\frac{n}{2}+\gamma$ with $\gamma\in (0,+\infty)\backslash \mathbb{N}$ and  denote the solution of \eqref{Sec 5 equ} by
		\begin{align*}
		u:=\mathcal{P}(s)f=\rho^{n-s}F+\rho^s G\;\; \;\textrm{with} \;\; F,G\in C^{\infty}(\overline{\mathbb{B}^{n+1}}).
		\end{align*}
		From our earlier analysis, we know that $F\big|_{\S^n}=f$, and we define the scattering operator as
		\begin{align*}
			S(s)f=G\big|_{\S^n}=\frac{\Gamma(-\gamma)}{2^{2\gamma}\Gamma(\gamma)}P_{2\gamma}f=\frac{1}{c_{\gamma}}P_{2\gamma}f,
		\end{align*}
		where $P_{2\gamma}$ is nothing but  the fractional GJMS operator of order $2\gamma$. Moreover, we have the following expansion
		\begin{align}\nonumber
			\mathcal{P}\left(\frac{n}{2}+\gamma\right)f
			=&\rho^{\frac{n}{2}-\gamma}\left(f+\rho^2f_1+\rho^4f_2+\cdots\right)+\\ \label{Expansion Ps}
&\frac{\Gamma(-\gamma)}{2^{2\gamma}\Gamma(\gamma)}\rho^{\frac{n}{2}+\gamma}
\left(P_{2\gamma}f+\rho^2g_1+\rho^4g_2+\cdots\right),
		\end{align}
with $f_i, g_j\in C^{\infty}(\mathbb{S}^{n})$, $i,j\in\mathbb{N}$.
	\end{rem}

	\subsection{Boundary operators}
	\subsubsection{The boundary operators with small indices}
For simplicity, we set
	\begin{align}\label{Sec 5 notation lap}
			\tilde{\Delta}_{+}:=\Delta_{+}+\frac{n^2}{4},
			\end{align}
			we also recall that the metric $g_{\B}$ can be rewrite as
			\begin{align}\label{Sec 5 notation g}
					g_{\B}=\frac{\ud \rho^{2}+(1-\frac{\rho^2}{4})^2g_{\S^n}}{\rho^2}, \quad\mathrm{where}\quad\rho=\frac{2(1-|x|)}{1+|x|}.
			\end{align}
Then the boundary operators with small indices is  as follows:
\begin{align*}	
B^{2\gamma}_{0}(U)&=U\big|_{\rho=0};\\
B^{2\gamma}_{2j}(U)&=\frac{1}{b_{2j}}\rho^{-\frac{n}{2}+\gamma-2j}\prod_{l=0}^{j-1}\left[
\left(\tilde{\Delta}_{+}-(\gamma-2l)^2\right)\left(\tilde{\Delta}_{+}-(\gamma+2l-2\lfloor \gamma\rfloor)^2\right)\right]\left(\rho^{\frac{n}{2}-\gamma}U\right)\Big|_{\rho=0},\\
&\qquad 1\leq j\leq \lfloor \gamma/2\rfloor;\\	B^{2\gamma}_{2j+2[\gamma]}(U)&=\frac{1}{b_{2j+2[\gamma]}}\rho^{-\frac{n}{2}+\gamma-2j-2[\gamma]}\prod_{l=0}^{j}\left(\tilde{\Delta}_{+}-(\gamma-2l)^2\right)\prod_{l=0}^{j-1}\left(\tilde{\Delta}_{+}-(\gamma+2l-2\lfloor \gamma\rfloor)^2\right)\left(\rho^{\frac{n}{2}-\gamma}U\right)\Big|_{\rho=0},\\
&\;\qquad 0\leq j\leq \lfloor \gamma\rfloor-\lfloor \gamma/2\rfloor-1,
\end{align*}
where $U\in\mathcal{C}^{2\gamma}(\overline{\B^{n+1}})$. For simplicity, we set
\begin{align*}
\tilde{B}^{2\gamma}_{2j}=&b_{2j}B^{2\gamma}_{2j},\qquad \qquad1\leq j\leq \lfloor \gamma/2\rfloor;\\
\tilde{B}^{2\gamma}_{2j+2[\gamma]}=&b_{2j+2[\gamma]}B^{2\gamma}_{2j+2[\gamma]},\qquad 0\leq j\leq \lfloor \gamma\rfloor-\lfloor \gamma/2\rfloor-1.
\end{align*}

To start, we present a simple and easily verifiable formula that illustrates how the operator $\tilde{\Delta}_{+}$ acts on expressions of the form $\rho^{\alpha}f$, where $f \in C^{\infty}(\S^n)$. This formula will be a valuable tool in our subsequent analysis.

\begin{lem}\label{Sec 5 Calculate lem 1}
	Let $f\in C^{\infty}(\S^n)$ and $\alpha\in \R$. Then  we have
	\begin{align*}
		\tilde{\Delta}_{+}(\rho^{\alpha}f)=\left(\alpha-\frac{n}{2}\right)^2\rho^{\alpha}f+O(\rho^{\alpha+2}).
	\end{align*}
	In particular, for $\gamma\in \R$ and $\alpha=\frac{n}{2}-\gamma$, it holds
	\begin{align*}
		\left(\tilde{\Delta}_{+}-\gamma^2\right)(\rho^{\frac{n}{2}-\gamma}f)=O(\rho^{\frac{n}{2}-\gamma+2}).
	\end{align*}
\end{lem}
\begin{pf}
From the \eqref{Sec 5 notation g}, we get $\sqrt{G}=\rho^{-(n+1)}\left(1-\rho^2/4\right)^n\sqrt{G_{\S^n}}$. By the definition of $\tilde{\Delta}_{+}$, we find
\begin{align}\label{Sec 5 Calculate lem 1 formula a}
	\tilde{\Delta}_{+}(\rho^{\alpha}f)	=&\frac{1}{\sqrt{G}}\partial_{\rho}\left(\rho^{2}\sqrt{G}\partial_{\rho}(\rho^{\alpha}f)\right)+\rho^{\alpha+2}
\left(1-\frac{\rho^2}{4}\right)^{-2}\Delta_{\S^n}f+\frac{n^2}{4}\rho^{\alpha}f\nonumber\\
	=&\alpha(\alpha+1)\rho^{\alpha}f-\alpha \rho^{\alpha+1}\left[\frac{n+1}{\rho}+\frac{n\rho/2}{1-\rho^2/4}\right]f+\rho^{\alpha+2}\left(1-\frac{\rho^2}{4}\right)^{-2}\Delta_{\S^n}f+\frac{n^2}{4}\rho^{\alpha}f\nonumber\\
	=&\left(\alpha-\frac{n}{2}\right)^2\rho^{\alpha}f+O(\rho^{\alpha+2}).
\end{align}
Note that the remaining terms in the final expression are the form $\rho^{\alpha+2m}$ for $m\in\mathbb{N}$.
\end{pf}

Building on Lemma \ref{Sec 5 Calculate lem 1}, we now present some important calculations regarding the operators $B^{2\gamma}_{2j}$ and $B^{2\gamma}_{2j+2[\gamma]}$, which will provide significant convenience for future computations.

\begin{lem}\label{Sec 5 Calculate lem 2}
	Let $f\in C^{\infty}(\S^n)$, for  $m\in\mathbb{N}$, there holds
	\begin{align*}
	B^{2\gamma}_{2j}(\rho^{2j}f)=f,\qquad B^{2\gamma}_{2j}(\rho^{2j+ 2m}f)=0\qquad\mathrm{for}\qquad 0\leq j\leq \lfloor \gamma/2\rfloor
	\end{align*}
	and
	\begin{align*}
	B^{2\gamma}_{2j+2[\gamma]}(\rho^{2j+2[\gamma]}f)=f,\qquad B^{2\gamma}_{2j+2[\gamma]}(\rho^{2j+2[\gamma]+ 2m}f)=0 \qquad\mathrm{for}\qquad 0\leq j\leq \lfloor \gamma\rfloor-\lfloor \gamma/2\rfloor-1.
	\end{align*}
\end{lem}
\begin{pf}
For any fixed $0\leq j\leq \lfloor \gamma/2\rfloor$,	using the Lemma \ref{Sec 5 Calculate lem 1}, we have
	\begin{align*}
			\left(\tilde{\Delta}_{+}-\gamma^2\right)(\rho^{\frac{n}{2}-\gamma+2j}f)=\left((\gamma-2j)^2-\gamma^2\right)
\rho^{\frac{n}{2}-\gamma+2j}f+O(\rho^{\frac{n}{2}-\gamma+2j+2}).
	\end{align*}
	Repeating the above computation, we obtain
	\begin{align}\label{Sec 5 Calculate lem 2 formula a}
		&\prod_{l=0}^{j-1}\left(\tilde{\Delta}_{+}-(\gamma-2l)^2\right)(\rho^{\frac{n}{2}-\gamma+2j}f)\nonumber\\
		=&\prod_{l=0}^{j-1}\left((\gamma-2j)^2-(\gamma-2l)^2\right)\rho^{\frac{n}{2}-\gamma+2j}f+O(\rho^{\frac{n}{2}-\gamma+2j+2}).
	\end{align}
	For the $\rho^{2j+2m}f$, we can also have
	\begin{align*}
		\rho^{-\frac{n}{2}+\gamma-2j}\prod_{l=0}^{j-1}\left(\tilde{\Delta}_{+}-(\gamma-2l)^2\right)\left(\tilde{\Delta}_{+}-(\gamma+2l-2\lfloor \gamma\rfloor)^2\right)(\rho^{2j+2m}f)=O(\rho^{2m}).
	\end{align*}
Therefore, we conclude that,  for $0\leq j\leq \lfloor \gamma/2\rfloor$, it holds
	\begin{align}\label{Sec 5 Calculate lem 2 formula b}
	B^{2\gamma}_{2j}(\rho^{2j}f)=\frac{1}{b_{2j}}	\tilde{B}^{2\gamma}_{2j}(\rho^{2j}f)=f \qquad \mathrm{and}\qquad
B^{2\gamma}_{2j}(\rho^{2j+ 2m}f)=\frac{1}{b_{2j}}\tilde{B}^{2\gamma}_{2j}(\rho^{2j+ 2m}f)=0.
	\end{align}
	By the similar calculation, for $0\leq j\leq \lfloor \gamma\rfloor-\lfloor \gamma/2\rfloor-1$, we can also obtain
	\begin{align}\label{Sec 5 Calculate lem 2 formula c}
	B^{2\gamma}_{2j+2[\gamma]}(\rho^{2j+2\lfloor \gamma\rfloor}f)=f \qquad \mathrm{and}\qquad 	B^{2\gamma}_{2j+2[\gamma]}(\rho^{2j+2\lfloor \gamma\rfloor+2m}f)=0.
	\end{align}
	This completes the proof of Lemma \ref{Sec 5 Calculate lem 2}.
\end{pf}

Next, we turn our attention to another type of calculation, specifically how the operator $B^{2\gamma}_{2j}$ acts on $\rho^{2[\gamma]+2m}f$ and how $B^{2\gamma}_{2j+2[\gamma]}$ acts on $\rho^{2m}f$.

\begin{lem}\label{Sec 5 Calculate lem 3}
	Suppose $f\in C^{\infty}(\S^n)$ and $m\in \mathbb{N}\cup\{0\}$, then we have \begin{align*}
			B^{2\gamma}_{2j}(\rho^{2[\gamma]+2m}f)=&	0, \qquad\mathrm{for}\qquad 0\leq j\leq \lfloor \gamma/2\rfloor;\\	
		B^{2\gamma}_{2j+2[\gamma]}(\rho^{2m}f)=&	0, \qquad\mathrm{for}\qquad 0\leq j\leq \lfloor \gamma\rfloor-\lfloor \gamma/2\rfloor-1.
	\end{align*}
\end{lem}
\begin{pf}
	It is enough to show  the case $m=0$ for both identities.

By the Lemma \ref{Sec 5 Calculate lem 1}, we have
	\begin{align*}
	\left(\tilde{\Delta}_{+}-(\gamma-2\lfloor \gamma\rfloor)^2\right)(\rho^{\frac{n}{2}-\gamma+2[\gamma]}f)=\rho^{\frac{n}{2}-\gamma+2[\gamma]+2}f_0+O(\rho^{\frac{n}{2}-\gamma+2[\gamma]+4})
	\end{align*}
	for some  function $f_0\in C^{\infty}(\S^n)$. Similarly, by using Lemma \ref{Sec 5 Calculate lem 1} again, we also get
	\begin{align*}
&\left(\tilde{\Delta}_{+}-(\gamma-2\lfloor \gamma\rfloor+2)^2\right)	\left(\tilde{\Delta}_{+}-(\gamma-2\lfloor \gamma\rfloor)^2\right)(\rho^{\frac{n}{2}-\gamma+2[\gamma]}f)\\
&=\rho^{\frac{n}{2}-\gamma+2[\gamma]+4}f_1+O(\rho^{\frac{n}{2}-\gamma+2[\gamma]+6})
	\end{align*}
	for some $f_1\in C^{\infty}(\S^n)$, hence it follows that, for some $f_2\in C^{\infty}(\S^n)$, it holds
	\begin{align*}
		\prod_{l=0}^{j-1}\left(\tilde{\Delta}_{+}-(\gamma+2l-2\lfloor \gamma\rfloor)^2\right)\left(\rho^{\frac{n}{2}-\gamma+2[\gamma]}f\right)=\rho^{\frac{n}{2}-\gamma+2[\gamma]+2j}f_2+O(\rho^{\frac{n}{2}-\gamma+2[\gamma]+2j+2 }).
	\end{align*}
	This will lead to
	\begin{align*}
\rho^{-\frac{n}{2}+\gamma-2j}\prod_{l=0}^{j}\left(\tilde{\Delta}_{+}-(\gamma-2l)^2\right)\prod_{l=0}^{j-1}\left(\tilde{\Delta}_{+}-(\gamma+2l-2\lfloor \gamma\rfloor)^2\right)(\rho^{2[\gamma]}f)=O(\rho^{2[\gamma]}),
	\end{align*}
	which implies  $B^{2\gamma}_{2j}(\rho^{2[\gamma]}f)=0$ for $0\leq j\leq \lfloor \gamma/2\rfloor$.

The proof of the second identity is exactly similar and we omit it.
This completes the proof of Lemma \ref{Sec 5 Calculate lem 3}.
	
\end{pf}

Based on the above computations, we can provide a clear description of the relationship between the expansion coefficients of a function in $\mathcal{C}^{2\gamma}(\overline{\B^{n+1}})$ and the boundary operators.
	
	\begin{cor}\label{Coro expansion}
		Suppose $f\in 		\mathcal{C}^{2\gamma}(\overline{\B^{n+1}})$, namely, it has expansion
		\begin{align*}
			f=\sum_{l=0}^{+\infty}\rho^{2l}f_{2l}+\sum_{l=0}^{+\infty}\rho^{2l+2[\gamma]}f_{2l+2[\gamma]},
		\end{align*}
		then
		\begin{align*}
				f_{2j}=B^{2\gamma}_{2j}\left(f-\sum_{l=0}^{j-1}\rho^{2l}f_{2l}\right)\qquad\mathrm{for}\qquad 0\leq j\leq \lfloor \gamma/2\rfloor;
		\end{align*}
		\begin{align*}
		f_{2j+2[\gamma]}=B^{2\gamma}_{2j+2[\gamma]}\left(f-\sum_{l=0}^{j-1}\rho^{2l+2[\gamma]}f_{2l+2[\gamma]}\right)\qquad\mathrm{for}\qquad 0\leq j\leq \lfloor \gamma\rfloor-\lfloor \gamma/2\rfloor-1.
		\end{align*}
	\end{cor}
	
	\subsubsection{Dirichlet problem}
	In this subsection, we begin by proving the existence and uniqueness of the polyharmonic Dirichlet problem. This allows us to define the polyharmonic extension for any function in $\mathcal{C}^{2\gamma}(\overline{\B^{n+1}})$, which is crucial for defining boundary operators with large indices.	To start, we recall the polyharmonic operators as follows: for $s=\frac{n}{2}+\gamma$, let
	\begin{align*}
	D_s=-\Delta_{+}-s(n-s),\qquad L_{2k}^{+}=\prod_{j=0}^{k-1}D_{s-2j}
	\end{align*}
	and
	\begin{align*}
		L_{2k}=\rho^{-\frac{n}{2}+\gamma-2k}\circ L^{+}_{2k}\circ \rho^{\frac{n}{2}-\gamma}
	\end{align*}
	
	\begin{thm}\label{Dirichlet Thm}
		Given $\gamma\in (0, +\infty)\backslash\mathbb{N}$ and boundary functions
		\begin{align*}
			f^{(2j)}\in C^{\infty}(\S^n) \qquad&\mathrm{for}\qquad 0\leq j\leq\lfloor \gamma/2\rfloor;\\
				\phi^{(2m)}\in C^{\infty}(\S^n) \qquad&\mathrm{for}\qquad 0\leq m\leq\lfloor \gamma\rfloor-\lfloor \gamma/2\rfloor-1.
		\end{align*}
		Then the Dirichlet problem $(k=\lfloor \gamma\rfloor+1)$
		\begin{align}\label{Dirichlet Thm equ}
			\begin{cases}
			\displaystyle L_{2k}U=0 \qquad&\mathrm{in}\qquad \B^{n+1},\\
			\displaystyle B^{2\gamma}_{2j}(U)=f^{(2j)}\qquad&\mathrm{for}\qquad 0\leq j\leq\lfloor \gamma/2\rfloor,\\
				\displaystyle B^{2\gamma}_{2m+2[\gamma]}(U)=\phi^{(2m)}\qquad&\mathrm{for}\qquad 0\leq m\leq\lfloor \gamma\rfloor-\lfloor \gamma/2\rfloor-1
			\end{cases}
		\end{align}
		has a unique solution
		\begin{align*}
		U=\sum_{j=0}^{\lfloor \gamma/2\rfloor}\rho^{-\frac{n}{2}+\gamma}\mathcal{P}\left(\frac{n}{2}+\gamma-2j\right)f^{(2j)}+\sum_{m=0}^{\lfloor \gamma\rfloor-\lfloor \gamma/2\rfloor-1}\rho^{-\frac{n}{2}+\gamma}\mathcal{P}\left(\frac{n}{2}+\lfloor \gamma\rfloor-[\gamma]-2m\right)\phi^{(2m)}.
		\end{align*}
	\end{thm}
	\begin{pf} The proof is exactly similar to that in \cite[Theorem 4.1]{Flynn&Lu&Yang1} and \cite[Theorem 4.1]{Flynn&Lu&Yang} and we omit it.
		\end{pf}
	
	\begin{definition}\label{Def polyharmonic}
	For any $ U \in \mathcal{C}^{2\gamma}(\overline{\B^{n+1}}) $, we say that $ \tilde{U} $ is a polyharmonic extension of $ U $ if $ \tilde{U} $ satisfies the following conditions
			\begin{align*}
		\begin{cases}
		\displaystyle L_{2k}\tilde{U}=0 \qquad&\mathrm{in}\qquad \B^{n+1},\\
		\displaystyle B^{2\gamma}_{2j}(\tilde{U})=B^{2\gamma}_{2j}(U)\qquad&\mathrm{for}\qquad 0\leq j\leq\lfloor \gamma/2\rfloor,\\
		\displaystyle B^{2\gamma}_{2m+2[\gamma]}(\tilde{U})=B^{2\gamma}_{2m+2[\gamma]}(U)\qquad&\mathrm{for}\qquad 0\leq m\leq\lfloor \gamma\rfloor-\lfloor \gamma/2\rfloor-1.
		\end{cases}
		\end{align*}
	\end{definition}

From Corollary \ref{Coro expansion}, we obtain the asymptotic expansion at the boundary
\begin{align}\label{U-U expansion}
	U-\tilde{U}\in \rho^{2(\lfloor \gamma/2\rfloor+1)}C^{\infty}_{\mathrm{even}}(\overline{\B^{n+1}})+\rho^{2(\gamma-\lfloor \gamma/2\rfloor)}C^{\infty}_{\mathrm{even}}(\overline{\B^{n+1}}).
\end{align}
 From Theorem \ref{Dirichlet Thm}, we can write
	
	\begin{align}\label{decom U}
		\tilde{U}=\sum_{j=0}^{\lfloor \gamma/2\rfloor}U_j+\sum_{j=0}^{\lfloor \gamma\rfloor-\lfloor \gamma/2\rfloor-1}U^{'}_j,
	\end{align}
	where
	\begin{align*}
		U_j=&\rho^{-\frac{n}{2}+\gamma}\mathcal{P}\left(\frac{n}{2}+\gamma-2j\right)B^{2\gamma}_{2j}U;\\
		U^{'}_j=&\rho^{-\frac{n}{2}+\gamma}\mathcal{P}\left(\frac{n}{2}+\lfloor \gamma\rfloor-[\gamma]-2j\right)B^{2\gamma}_{2j+2[\gamma]}U.
	\end{align*}
	For $j=0,1,\cdots,\lfloor \gamma\rfloor-\lfloor \gamma/2\rfloor-1$, define $U_{\lfloor \gamma/2\rfloor+j+1}=U^{'}_{\lfloor \gamma\rfloor-\lfloor \gamma/2\rfloor-1-j}$. We get
	\begin{align*}
		U_{\lfloor \gamma/2\rfloor+j+1}=\rho^{-\frac{n}{2}+\gamma}\mathcal{P}\left(\frac{n}{2}-\gamma+2(\lfloor \gamma/2\rfloor+j+1)\right)B^{2\gamma}_{2\gamma-2(\lfloor \gamma/2\rfloor+j+1)}U.
	\end{align*}
From the above analysis, it is clear that
	\begin{align*}
	\left(\tilde{\Delta}_{+}-(\gamma-2j)^2\right)\left(\rho^{\frac{n}{2}-\gamma}U_j\right)=0\qquad\mathrm{for}\qquad 0\leq j\leq \lfloor \gamma\rfloor.
	\end{align*}
	Using Remark \ref{Vital Rem}, we obtain
	\begin{align}\label{decom U1}
		U_j=\begin{cases}
		\displaystyle \rho^{2j}F_j+\rho^{2\gamma-2j}G_j, \qquad j=0,1,\cdots, \lfloor \gamma/2\rfloor,\\
		\displaystyle \rho^{2\gamma-2j}F_j+\rho^{2j}G_j, \qquad j=\lfloor \gamma/2\rfloor+1,\cdots,\lfloor \gamma\rfloor,
		\end{cases}
	\end{align}
	where $F_j, G_j\in C^{\infty}(\overline{\B^{n+1}})$,
	\begin{align}\label{asy Fj}
		F_j\big|_{\rho=0}=\begin{cases}
		\displaystyle B^{2\gamma}_{2j}U, \qquad\quad~ j=0,1,\cdots, \lfloor \gamma/2\rfloor,\\
		\displaystyle B^{2\gamma}_{2\gamma-2j}U, \qquad j=\lfloor \gamma/2\rfloor+1,\cdots,\lfloor \gamma\rfloor
		\end{cases}
	\end{align}
	and
	\begin{align*}
			G_j\big|_{\rho=0}=\begin{cases}
		\displaystyle S\left(n/2+\gamma-2j\right)B^{2\gamma}_{2j}U, \qquad\quad~ j=0,1,\cdots, \lfloor \gamma/2\rfloor,\\
		\displaystyle S\left(n/2-\gamma+2j\right)B^{2\gamma}_{2\gamma-2j}U, \qquad j=\lfloor \gamma/2\rfloor+1,\cdots,\lfloor \gamma\rfloor.
		\end{cases}
	\end{align*}

	\subsubsection{The boundary operators with large indices}\label{sec large index}
	In this subsection, we begin by defining the boundary operators with large indices. For notational convenience, we introduce the following definitions:
	\begin{align*}
		\Pi^{\gamma}=(-1)^{\lfloor \gamma\rfloor+1}\prod_{l=0}^{\lfloor \gamma\rfloor}\left(\tilde{\Delta}_{+}-(\gamma-2l)^2\right)=L^{+}_{2k}, \qquad k=\lfloor \gamma\rfloor+1
	\end{align*}
	and
	\begin{align*}
		\Pi^{\gamma}_{j}=&(-1)^{\lfloor \gamma\rfloor}\prod_{l=0}^{j-1}\left(\tilde{\Delta}_{+}-(\gamma-2l)^2\right)\prod_{l=j+1}^{\lfloor \gamma\rfloor}\left(\tilde{\Delta}_{+}-(\gamma-2l)^2\right)\\
		=&(-1)^{\lfloor \gamma\rfloor}\prod_{l=\lfloor \gamma\rfloor-j+1}^{\lfloor \gamma\rfloor}\left(\tilde{\Delta}_{+}-(\gamma-2l-2[\gamma])^2\right)\prod_{l=0}^{\lfloor \gamma\rfloor-j-1}\left(\tilde{\Delta}_{+}-(\gamma-2l-2[\gamma])^2\right),
	\end{align*}
which can also be written as
	\begin{align*}
		\Pi^{\gamma}_j=\frac{\Pi^{\gamma}}{-\left(\tilde{\Delta}_{+}-(\gamma-2j)^2\right)}.
	\end{align*}
	For the operator $ \Pi^{\gamma}_j $, we list the following properties, which can be derived through similar computations as those in Lemma \ref{Sec 5 Calculate lem 2} and Lemma \ref{Sec 5 Calculate lem 3}.
	
	\begin{lem}\label{Sec 5 Pi act formula Lem}
		For $0\leq j\leq \lfloor \gamma\rfloor$, if $f\in \rho^{2i}C^{\infty}_{\mathrm{even}}(\overline{\B^{n+1}})$, then
		\begin{align*}
			\Pi^{\gamma}_j\left(\rho^{\frac{n}{2}-\gamma}f\right)\in& \rho^{\frac{n}{2}-\gamma+2j}C^{\infty}_{\mathrm{even}}(\overline{\B^{n+1}})\quad~~~ \qquad\mathrm{for}\qquad i\leq j,\\
			\Pi^{\gamma}_j\left(\rho^{\frac{n}{2}-\gamma}f\right)\in& \rho^{\frac{n}{2}-\gamma+2\lfloor \gamma\rfloor+2}C^{\infty}_{\mathrm{even}}(\overline{\B^{n+1}}) \qquad\mathrm{for}\qquad i> j.
		\end{align*}
		If $f\in \rho^{2i+2[\gamma]}C^{\infty}_{\mathrm{even}}(\overline{\B^{n+1}})$, then
		\begin{align*}
				\Pi^{\gamma}_j\left(\rho^{\frac{n}{2}-\gamma}f\right)\in& \rho^{\frac{n}{2}+\gamma-2j}C^{\infty}_{\mathrm{even}}(\overline{\B^{n+1}})\quad~~ \qquad\mathrm{for}\qquad i\leq \lfloor \gamma\rfloor-j,\\
					\Pi^{\gamma}_j\left(\rho^{\frac{n}{2}-\gamma}f\right)\in& \rho^{\frac{n}{2}+\gamma+2}C^{\infty}_{\mathrm{even}}(\overline{\B^{n+1}})\quad~~~ \qquad\mathrm{for}\qquad i>\lfloor \gamma\rfloor-j.
		\end{align*}
	\end{lem}
Clearly, combining \eqref{U-U expansion} with Lemma \ref{Sec 5 Pi act formula Lem}, we have
\begin{itemize}
\item $0\leq j\leq \lfloor \gamma/2\rfloor$,
\begin{align}\label{Pi asy 1}
	\Pi^{\gamma}_j\left(\rho^{\frac{n}{2}-\gamma}(U-\tilde{U})\right)\in \rho^{\frac{n}{2}-\gamma+2\lfloor \gamma\rfloor+2}C^{\infty}_{\mathrm{even}}(\overline{\B^{n+1}})+\rho^{\frac{n}{2}+\gamma-2j}C^{\infty}_{\mathrm{even}}(\overline{\B^{n+1}}).
\end{align}
\item $\lfloor \gamma/2\rfloor+1\leq j\leq \lfloor \gamma\rfloor$,
\begin{align}\label{Pi asy 2}
\Pi^{\gamma}_j\left(\rho^{\frac{n}{2}-\gamma}(U-\tilde{U})\right)\in \rho^{\frac{n}{2}-\gamma+2j}C^{\infty}_{\mathrm{even}}(\overline{\B^{n+1}})+\rho^{\frac{n}{2}+\gamma+2}C^{\infty}_{\mathrm{even}}(\overline{\B^{n+1}}).
\end{align}
\end{itemize}
	For any $ U \in \mathcal{C}^{2\gamma}(\overline{\B^{n+1}}) $, let $ \tilde{U} $ denote its polyharmonic extension. Following \cite{Flynn&Lu&Yang1}, we define the boundary operators with large indices as follows:
\begin{itemize}
	\item $1+\lfloor \gamma/2\rfloor\leq j\leq \lfloor \gamma\rfloor$,
	\begin{align*}		B^{2\gamma}_{2j}(U)=\frac{1}{b_{2j}}\rho^{-\frac{n}{2}+\gamma-2j}\Pi^{\gamma}_j\left(\rho^{\frac{n}{2}-\gamma}(U-\tilde{U})\right)\Big|_{\rho=0}
+\frac{1}{c_{2j-\gamma}}P_{4j-2\gamma}B^{2\gamma}_{2\gamma-2j}(U).
	\end{align*}
	\item  $\lfloor \gamma\rfloor-\lfloor \gamma/2\rfloor\leq j\leq \lfloor \gamma\rfloor$,
	\begin{align*}
	B^{2\gamma}_{2j+2[\gamma]}(U)=&\frac{1}{b_{2j+2[\gamma]}}\rho^{-\frac{n}{2}+\gamma-2j-2[\gamma]}\Pi^{\gamma}_{\lfloor \gamma\rfloor-j}\left(\rho^{\frac{n}{2}-\gamma}(U-\tilde{U})\right)\Big|_{\rho=0}\\
	&+\frac{1}{c_{2j+[\gamma]-\lfloor \gamma\rfloor}}P_{4j+2[\gamma]-2\lfloor \gamma\rfloor}B^{2\gamma}_{2\lfloor \gamma\rfloor-2j}(U),
	\end{align*}
\end{itemize}	
where
\begin{align*}
	b_{2j}=&(-1)^{\lfloor \gamma\rfloor}\prod_{i=0,i\not=j}^{\lfloor \gamma\rfloor}\left((\gamma-2j)^2-(\gamma-2i)^2\right),\\
	b_{2j+2[\gamma]}=&(-1)^{\lfloor \gamma\rfloor}\prod_{i=0,i\not=\lfloor \gamma\rfloor-j}^{\lfloor \gamma\rfloor}\left((\gamma-2j-2[\gamma])^2-(\gamma-2i)^2\right)
\end{align*}
are chosen such that
\begin{align*}
	B^{2\gamma}_{2j}(\rho^{2j})= B^{2\gamma}_{2j+2[\gamma]}(\rho^{2j+2[\gamma]})=1.
\end{align*}
We can also express the boundary operators $B^{2\gamma}_{2j+2[\gamma]}$  in the following form
\begin{align}\label{Sec 5 formula B 2-2}
	B^{2\gamma}_{2\gamma-2j}(U)=\frac{1}{b_{2\gamma-2j}}\rho^{-\frac{n}{2}-\gamma+2j}\Pi^{\gamma}_{j}
\left(\rho^{\frac{n}{2}-\gamma}(U-\tilde{U})\right)\Big|_{\rho=0}+\frac{1}{c_{\gamma-2j}}P_{2\gamma-4j}B^{2\gamma}_{2j}(U),
\end{align}
where $0\leq j\leq \lfloor \gamma/2\rfloor$.

	\subsection{Proof of theorem \ref{Conformal Thm}  and  \ref{Intrinsic Thm}}
	
	\textbf{Proof of Theorem \ref{Conformal Thm}:} The proof is similar to that given in \cite[Theorem 1.4]{Flynn&Lu&Yang1}.
For readers convenience, we describe it as follows.

Firstly, the conformal invariance of boundary operators with small indices is straightforward. From the definitions of
	 $	\widehat{B}^{2\gamma}_{2j}(U)$ and $\widehat{B}^{2\gamma}_{2j+2[\gamma]}(U)$, we can see that
	\begin{align*}
			\widehat{B}^{2\gamma}_{2j}(U)&=\frac{1}{b_{2j}}\widehat{\tilde{B}}^{2\gamma}_{2j}(U) \qquad\mathrm{for}\qquad 0\leq j\leq \lfloor \gamma/2\rfloor;\\
			\widehat{B}^{2\gamma}_{2j+2[\gamma]}(U)&=\frac{1}{b_{2j+2[\gamma]}}\widehat{\tilde{B}}^{2\gamma}_{2j+2[\gamma]}(U) \qquad\mathrm{for}\qquad 0\leq j\leq \lfloor \gamma\rfloor-\lfloor \gamma/2\rfloor-1,
	\end{align*}
	where
	\begin{align*}
	\widehat{\tilde{B}}^{2\gamma}_{2j}(U)&=\hat{\rho}^{-\frac{n}{2}+\gamma-2j}\prod_{l=0}^{j-1}\left(\tilde{\Delta}_{+}-(\gamma-2l)^2\right)
\left(\tilde{\Delta}_{+}-(\gamma+2l-2\lfloor \gamma\rfloor)^2\right)\left(\hat{\rho}^{\frac{n}{2}-\gamma}U\right)\Big|_{\rho=0};\\
	\widehat{\tilde{B}}^{2\gamma}_{2j+2[\gamma]}(U)&=\hat{\rho}^{-\frac{n}{2}+\gamma-2j-2[\gamma]}\prod_{l=0}^{j}\left(\tilde{\Delta}_{+}-(\gamma-2l)^2\right)\prod_{l=0}^{j-1}\left(\tilde{\Delta}_{+}-(\gamma+2l-2\lfloor \gamma\rfloor)^2\right)\left(\hat{\rho}^{\frac{n}{2}-\gamma}U\right)\Big|_{\rho=0}.
	\end{align*}
	Hence, the conformally invariant identities \eqref{Conformal Thm property 1} and \eqref{Conformal Thm property 2} follow from $\hat{\rho}=e^{\tau}\rho$.
	
	Secondly, for large $j$, we only need to address
	 $B^{2\gamma}_{2j}$,
	 the remaining part can be obtained through a similar argument. Now, from the definition of $\widehat{B}^{2\gamma}_{2j}$,
	\begin{align}\label{Sec 5 hat B}
		\widehat{B}^{2\gamma}_{2j}(U)=&\frac{1}{b_{2j}}\hat{\rho}^{-\frac{n}{2}+\gamma-2j}\Pi^{\gamma}_j\left(\hat{\rho}^{\frac{n}{2}-\gamma}(U-\widehat{\tilde{U}})\right)
\Big|_{\rho=0}+\frac{1}{c_{2j-\gamma}}\widehat{P}_{4j-2\gamma}\widehat{B}^{2\gamma}_{2\gamma-2j}(U),
	\end{align}
	where $\widehat{\tilde{U}}$ is the unique solution of
		\begin{align*}
	\begin{cases}
	\displaystyle \widehat{L}_{2k}\widehat{\tilde{U}}=0 \qquad&\mathrm{in}\qquad \B^{n+1},\\
	\displaystyle \widehat{B}^{2\gamma}_{2j}(\widehat{\tilde{U}})=\widehat{B}^{2\gamma}_{2j}(U)\qquad&\mathrm{for}\qquad 0\leq j\leq\lfloor \gamma/2\rfloor,\\
	\displaystyle \widehat{B}^{2\gamma}_{2m+2[\gamma]}(\widehat{\tilde{U}})=\widehat{B}^{2\gamma}_{2m+2[\gamma]}(U)\qquad&\mathrm{for}\qquad 0\leq m\leq\lfloor \gamma\rfloor-\lfloor \gamma/2\rfloor-1.
	\end{cases}
	\end{align*}
One crucial observation that significantly aids our proof is the following claim:	
	\begin{align}\label{Sec 5 claim 1}
		\widehat{\tilde{U}}=e^{-\left(\frac{n}{2}-\gamma\right)\tau}(e^{\left(\frac{n}{2}-\gamma\right)\tau}U)^{\widetilde{}},
	\end{align}
	where $(e^{\left(\frac{n}{2}-\gamma\right)\tau}U)^{\widetilde{}}$ is the solution of
		\begin{align*}
	\begin{cases}
	\displaystyle L_{2k}(e^{\left(\frac{n}{2}-\gamma\right)\tau}U)^{\widetilde{}}=0 \qquad&\mathrm{in}\qquad \B^{n+1},\\
	\displaystyle B^{2\gamma}_{2j}((e^{\left(\frac{n}{2}-\gamma\right)\tau}U)^{\widetilde{}})=B^{2\gamma}_{2j}(e^{\left(\frac{n}{2}-\gamma\right)\tau}U)\qquad&\mathrm{for}\qquad 0\leq j\leq\lfloor \gamma/2\rfloor,\\
	\displaystyle B^{2\gamma}_{2m+2[\gamma]}((e^{\left(\frac{n}{2}-\gamma\right)\tau}U)^{\widetilde{}})=B^{2\gamma}_{2m+2[\gamma]}(e^{\left(\frac{n}{2}-\gamma\right)\tau}U)\qquad&\mathrm{for}\qquad 0\leq m\leq\lfloor \gamma\rfloor-\lfloor \gamma/2\rfloor-1.
	\end{cases}
	\end{align*}
To prove \eqref{Sec 5 claim 1}, we first verify whether it satisfies the interior equation
	\begin{align*}
		\widehat{L}_{2k}\left(e^{-\left(\frac{n}{2}-\gamma\right)\tau}(e^{\left(\frac{n}{2}-\gamma\right)\tau}U)^{\widetilde{}}\right)=&\widehat{\rho}^{\;-\frac{n}{2}+\gamma-2k}\circ L^{+}_{2k} \left(\widehat{\rho}^{\frac{n}{2}-\gamma}e^{-\left(\frac{n}{2}-\gamma\right)\tau}(e^{\left(\frac{n}{2}-\gamma\right)\tau}U)^{\widetilde{}}\right)\\
		=&e^{\left(-\frac{n}{2}+\gamma-2k\right)\tau}L_{2k}\left((e^{\left(\frac{n}{2}-\gamma\right)\tau}U)^{\widetilde{}}\right)=0.
	\end{align*}
	Secondly, the boundary conditions can be verified through the following computations,
	\begin{align*}
		\widehat{B}^{2\gamma}_{2j}(e^{-\left(\frac{n}{2}-\gamma\right)\tau}(e^{\left(\frac{n}{2}-\gamma\right)\tau}U)^{\widetilde{}})=&e^{\left(-\frac{n}{2}+\gamma-2j\right)\tau}B^{2\gamma}_{2j}\left((e^{\left(\frac{n}{2}-\gamma\right)\tau}U)^{\widetilde{}}\right)\\
		=&e^{\left(-\frac{n}{2}+\gamma-2j\right)\tau}B^{2\gamma}_{2j}\left(e^{\left(\frac{n}{2}-\gamma\right)\tau}U\right)=\widehat{B}^{2\gamma}_{2j}(U)
	\end{align*}
	and
	\begin{align*}
		\widehat{B}^{2\gamma}_{2m+2[\gamma]}(e^{-\left(\frac{n}{2}-\gamma\right)\tau}(e^{\left(\frac{n}{2}-\gamma\right)\tau}U)^{\widetilde{}})=&e^{\left(-\frac{n}{2}+\gamma-2m-2[\gamma]\right)\tau}B^{2\gamma}_{2m+2[\gamma]}\left((e^{\left(\frac{n}{2}-\gamma\right)\tau}U)^{\widetilde{}}\right)\\
		=&e^{\left(-\frac{n}{2}+\gamma-2m-2[\gamma]\right)\tau}B^{2\gamma}_{2m+2[\gamma]}(e^{\left(\frac{n}{2}-\gamma\right)\tau}U)\\
		=&\widehat{B}^{2\gamma}_{2m+2[\gamma]}(U).
	\end{align*}
By the uniqueness result in Theorem \ref{Dirichlet Thm}, we complete the proof of claim \eqref{Sec 5 claim 1}. Now, returning to formula \eqref{Sec 5 hat B}, we have
	\begin{align}\label{Sec 5 thm 1 formula B}
			\widehat{B}^{2\gamma}_{2j}(U)=&\frac{1}{b_{2j}}e^{\left(-\frac{n}{2}+\gamma-2j\right)\tau}\rho^{-\frac{n}{2}+\gamma-2j}\Pi^{\gamma}_j\left(\rho^{\frac{n}{2}-\gamma}\left(e^{\left(\frac{n}{2}-\gamma\right)\tau}U-(e^{\left(\frac{n}{2}-\gamma\right)\tau}U)^{\widetilde{}}\right)\right)\Big|_{\rho=0}\nonumber
			\\
			&+\frac{1}{c_{2j-\gamma}}\widehat{P}_{4j-2\gamma}\left[e^{\left(-\frac{n}{2}-\gamma+2j\right)\tau|_{\S^n}}
B^{2\gamma}_{2\gamma-2j}(e^{\left(\frac{n}{2}-\gamma\right)\tau}U)\right].
	\end{align}
	Using the conformal property of GJMS operator $P_{2\gamma}$, we have
	\begin{align}\label{Sec 5 GJMS conformal}
		\widehat{P}_{4j-2\gamma}(f)=e^{\left(-\frac{n}{2}+\gamma-2j\right)\tau|_{\S^n}}P_{4j-2\gamma}\left(e^{\left(\frac{n}{2}+\gamma-2j\right)\tau|_{\S^n}}f\right) \qquad\mathrm{for}\qquad f\in C^{\infty}(\S^n).
	\end{align}
	Plugging \eqref{Sec 5 GJMS conformal} into \eqref{Sec 5 thm 1 formula B}, we conclude that
	\begin{align*}
			\widehat{B}^{2\gamma}_{2j}(U)=e^{\left(-\frac{n}{2}+\gamma-2j\right)\tau|_{\S^n}}B^{2\gamma}_{2j}\left(e^{\left(\frac{n}{2}-\gamma\right)\tau}U\right).
	\end{align*}
	\\

	\textbf{Proof of Theorem \ref{Intrinsic Thm}:}  Noticing that if $L_{2k}U=0$, then $U=\tilde{U}$, hence we have
	\begin{align*}
		B^{2\gamma}_{2\gamma-2j}(U)=\frac{1}{c_{\gamma-2j}}P_{{2\gamma-4j}}B^{2\gamma}_{2j}(U)
	\end{align*}
	and
	\begin{align*}
		B^{2\gamma}_{2\lfloor \gamma\rfloor-2j}(U)=\frac{1}{c_{\lfloor \gamma\rfloor-[\gamma]-2j}}P_{2\lfloor \gamma\rfloor-2[\gamma]-4j}B^{2\gamma}_{2[\gamma]+2j}(U).
	\end{align*}

	\subsection{Proof of theorem \ref{Symmetric Thm} and  \ref {Functional inequ Thm} }
	In this subsection, we will establish an important integral identity (Lemma \ref{Vital integral identity Lem}) and demonstrate the symmetry of the Dirichlet form (Theorem \ref{Symmetric Thm}). Building on this, we will complete the proofs of the remaining theorems. We begin with a simple integral identity, which can be viewed as a Green-type formula on the Poincar\'e ball.
	\begin{lem}\label{Green identity}
		Suppose $U, V\in 	\mathcal{C}^{2\gamma}(\overline{\B^{n+1}})$, then we have
		\begin{align*}
		\int_{\B^{n+1}}(U\tilde{\Delta}_{+}V-V\tilde{\Delta}_{+}U)\ud V_{g_{\B}}=-\int_{\S^n}\left(\rho^{1-n}U\partial_{\rho}V-\rho^{1-n}V\partial_{\rho}U\right)\big|_{\rho=0}\ud V_{\S^n}.
		\end{align*}
	\end{lem}
	\begin{pf}
		Set $\sqrt{G}:=\sqrt{G_{g_{\B}}}=\rho^{-(n+1)}\sqrt{G_g}$, then we can get
		\begin{align*}
				&\int_{\B^{n+1}}U\tilde{\Delta}_{+}V-V\tilde{\Delta}_{+}U\ud V_{g_{\B}}\\
				=&\int_{\S^n}\int_{0}^{2}\frac{1}{\sqrt{G_{\mathbb{S}^n}}}\left(U\partial_{\rho}(\rho^{1-n}\sqrt{G_{g}}\partial_{\rho}V)
-V\partial_{\rho}(\rho^{1-n}\sqrt{G_{g}}\partial_{\rho}U)\right)\ud\rho\ud V_{\S^n}\\
				=&-\int_{\S^n}\left(\rho^{1-n}U\partial_{\rho}V-\rho^{1-n}V\partial_{\rho}U\right)\big|_{\rho=0}\ud V_{\S^n},
		\end{align*}
		here we use $\sqrt{G_{g}}=\left(1-\frac{\rho^2}{4}\right)^n\sqrt{G_{\S^n}}$ and $\sqrt{G_{g}}\to 0$ as $\rho\to 2$.
	\end{pf}

By applying the Green-type formula and performing a careful asymptotic analysis near the boundary, we can establish the following nontrivial integral identity.
	
	\begin{lem}\label{Vital integral identity Lem}
		Suppose $U, V\in 	\mathcal{C}^{2\gamma}(\overline{\B^{n+1}})$ and $k=\lfloor \gamma\rfloor+1$, there holds
		\begin{align*}
			&\int_{\B^{n+1}}\rho^{1-2[\gamma]}UL_{2k}V\ud V_{g}\\
=&
\int_{\B^{n+1}}	\rho^{\frac{n}{2}-\gamma}(U-\tilde{U})L_{2k}^{+}\left(\rho^{\frac{n}{2}-\gamma}\left(V-\tilde{V}\right)\right)\ud V_{g_{\B}}\\
+&\sum_{j=0}^{\lfloor \gamma/2\rfloor}\sigma_{j,\gamma}\int_{\S^n}B^{2\gamma}_{2j}UB^{2\gamma}_{2\gamma-2j}V\ud V_{\S^n}+\sum_{j=\lfloor \gamma/2\rfloor+1}^{\lfloor \gamma\rfloor}\sigma_{j,\gamma}\int_{\S^n}B^{2\gamma}_{2j}VB^{2\gamma}_{2\gamma-2j}U\ud V_{\S^n}\\
+&\sum_{j=0}^{\lfloor \gamma/2\rfloor}\upzeta_{j,\gamma}\int_{\S^n}B^{2\gamma}_{2j}UP_{2\gamma-4j}B^{2\gamma}_{2j}V\ud V_{\S^n}+\sum_{j=\lfloor \gamma/2\rfloor+1}^{\lfloor \gamma\rfloor}\upzeta_{j,\gamma}\int_{\S^n}B^{2\gamma}_{2\gamma-2j}UP_{4j-2\gamma}B^{2\gamma}_{2\gamma-2j}V\ud V_{\S^n},
	\end{align*}
	where
	\begin{align*}
		\sigma_{j,\gamma}=&\begin{cases}
		\displaystyle 2^{2\lfloor \gamma\rfloor+1}j!(\lfloor \gamma\rfloor-j)!\frac{\Gamma(\gamma+1-j)}{\Gamma(\gamma-2j)}\frac{\Gamma(j+1-[\gamma])}{\Gamma(2j+1-\gamma)}, \quad 0\leq j\leq \lfloor \gamma/2\rfloor\\
		\displaystyle -2^{2\lfloor \gamma\rfloor+1}j!(\lfloor \gamma\rfloor-j)!\frac{\Gamma(\gamma+1-j)}{\Gamma(\gamma-2j)}\frac{\Gamma(j+1-[\gamma])}{\Gamma(2j+1-\gamma)},\quad \lfloor \gamma/2\rfloor+1\leq j\leq \lfloor \gamma\rfloor,
		\end{cases}\\
		\upzeta_{j,\gamma}=&\begin{cases}
		\displaystyle 2^{4j-2[\gamma]+1}j!(\lfloor \gamma\rfloor-j)!\frac{\Gamma(\gamma+1-j)}{\Gamma(\gamma+1-2j)}\frac{\Gamma(j+1-[\gamma])}{\Gamma(\gamma-2j)}, \quad 0\leq j\leq \lfloor \gamma/2\rfloor\\
		\displaystyle 2^{4\gamma-2[\gamma]+1}j!(\lfloor \gamma\rfloor-j)!\frac{\Gamma(\gamma+1-j)}{\Gamma(2j-\gamma)}\frac{\Gamma(j+1-[\gamma])}{\Gamma(2j+1-\gamma)},\quad \lfloor \gamma/2\rfloor+1\leq j\leq \lfloor \gamma\rfloor
		\end{cases}
	\end{align*}
	and $\upzeta_{j,\gamma}>0$ for all $0\leq j\leq \lfloor \gamma\rfloor$.
	\end{lem}
	\begin{pf}
		By the definition of $L_{2k} $ and $ L^{+}_{2k} \left( \rho^{\frac{n}{2}-\gamma} \tilde{V} \right) = 0 $, we can rewrite the following expression
		\begin{align*}
			\int_{\B^{n+1}}\rho^{1-2[\gamma]}UL_{2k}V\ud V_{g}=&\int_{\B^{n+1}}\rho^{\frac{n}{2}-\gamma}UL^{+}_{2k}(\rho^{\frac{n}{2}-\gamma
		}V)\ud V_{g_{\B}}\\
			=&\int_{\B^{n+1}}\rho^{\frac{n}{2}-\gamma}(U-\tilde{U})L^{+}_{2k}(\rho^{\frac{n}{2}-\gamma}(V-\tilde{V}))\ud V_{g_{\B}}\\
			+&\int_{\B^{n+1}}\rho^{\frac{n}{2}-\gamma}\tilde{U}L^{+}_{2k}(\rho^{\frac{n}{2}-\gamma}(V-\tilde{V}))\ud V_{g_{\B}}\\
			=:&A_1+A_2.
		\end{align*}
	Recall the decomposition of $ \tilde{U} $, as shown in equations \eqref{decom U} and \eqref{decom U1}, i.e,
		\begin{align}\label{Sec 5 Uj}
			\tilde{U}=\sum_{j=0}^{\lfloor \gamma\rfloor}U_j \qquad \mathrm{and}\qquad U_j=\begin{cases}
			\displaystyle \rho^{2j}F_j+\rho^{2\gamma-2j}G_j \qquad j=0,1,\cdots, \lfloor \gamma/2\rfloor,\\
			\displaystyle \rho^{2\gamma-2j}F_j+\rho^{2j}G_j\qquad j=\lfloor \gamma/2\rfloor+1,\cdots,\lfloor \gamma\rfloor.
			\end{cases}
		\end{align}
		Using Lemma \ref{Green identity} and the fact that $ \left( \tilde{\Delta}_{+} - (\gamma - 2j)^2 \right) \left( \rho^{\frac{n}{2} - \gamma} U_j \right) = 0 $, we proceed with
		\begin{align*}
		&\int_{\B^{n+1}}\rho^{\frac{n}{2}-\gamma}U_jL^{+}_{2k}(\rho^{\frac{n}{2}-\gamma}(V-\tilde{V}))\ud V_{g_{\B}}\\
			=&-\int_{\B^{n+1}}\rho^{\frac{n}{2}-\gamma}U_j\left(\tilde{\Delta}_{+}-(\gamma-2j)^2\right)\Pi^{\gamma}_{j}(\rho^{\frac{n}{2}-\gamma}(V-\tilde{V}))\ud V_{g_{\B}}\\
			=&\int_{\B^{n+1}}\left(\tilde{\Delta}_{+}-(\gamma-2j)^2\right)(\rho^{\frac{n}{2}-\gamma}U_j)\Pi^{\gamma}_{j}(\rho^{\frac{n}{2}-\gamma}(V-\tilde{V}))\ud V_{g_{\B}}\\
			-&\int_{\B^{n+1}}\rho^{\frac{n}{2}-\gamma}U_j\left(\tilde{\Delta}_{+}-(\gamma-2j)^2\right)\Pi^{\gamma}_{j}(\rho^{\frac{n}{2}-\gamma}(V-\tilde{V}))\ud V_{g_{\B}}\\
			=&-\int_{\S^n}\rho^{1-n}\partial_{\rho}(\rho^{\frac{n}{2}-\gamma}U_j)\Pi^{\gamma}_{j}(\rho^{\frac{n}{2}-\gamma}(V-\tilde{V}))\Big|_{\rho=0}\ud V_{\S^n}\\
			+&\int_{\S^n}\rho^{1-\frac{n}{2}-\gamma}U_j\partial_{\rho}\left[\Pi^{\gamma}_{j}(\rho^{\frac{n}{2}-\gamma}(V-\tilde{V}))\right]\Big|_{\rho=0}\ud V_{\S^n}.
		\end{align*}
		Therefore, we can see
		\begin{align*}
			A_2=&\sum_{j=0}^{\lfloor \gamma\rfloor}\int_{\S^n}\rho^{1-\frac{n}{2}-\gamma}U_j\partial_{\rho}\left[\Pi^{\gamma}_{j}(\rho^{\frac{n}{2}-\gamma}(V-\tilde{V}))\right]\Big|_{\rho=0}\ud V_{\S^n}\\
			-&\sum_{j=0}^{\lfloor \gamma\rfloor}\int_{\S^n}\rho^{1-n}\partial_{\rho}(\rho^{\frac{n}{2}-\gamma}U_j)\Pi^{\gamma}_{j}(\rho^{\frac{n}{2}-\gamma}(V-\tilde{V}))\Big|_{\rho=0}\ud V_{\S^n}.
		\end{align*}
		From equations \eqref{asy Fj}, \eqref{Pi asy 1}, (\ref{Sec 5 formula B 2-2}) and \eqref{Sec 5 Uj}, for $ 0 \leq j \leq \lfloor \gamma/2 \rfloor $, we can show
		\begin{align*}
			&\rho^{1-n}\partial_{\rho}(\rho^{\frac{n}{2}-\gamma}U_j)\Pi^{\gamma}_{j}(\rho^{\frac{n}{2}-\gamma}(V-\tilde{V}))\big|_{\rho=0}\\
		=&\left(\frac{n}{2}-\gamma+2j\right)F_j\big|_{\rho=0}\cdot\rho^{-\frac{n}{2}-\gamma+2j}\Pi^{\gamma}_{j}(\rho^{\frac{n}{2}-\gamma}(V-\tilde{V}))\big|_{\rho=0}\\
			&+\left(\frac{n}{2}+\gamma-2j\right)G_j\big|_{\rho=0}\cdot\rho^{-\frac{n}{2}+\gamma-2j}\Pi^{\gamma}_{j}(\rho^{\frac{n}{2}-\gamma}(V-\tilde{V}))\big|_{\rho=0}\\
			=&\left(\frac{n}{2}-\gamma+2j\right)B^{2\gamma}_{2j}(U)\cdot\rho^{-\frac{n}{2}+\gamma-2j}\Pi^{\gamma}_{j}(\rho^{\frac{n}{2}-\gamma}(V-\tilde{V}))\big|_{\rho=0}\\
			=&\left(\frac{n}{2}-\gamma+2j\right)b_{2\gamma-2j}B^{2\gamma}_{2j}(U)\cdot\left(B^{2\gamma}_{2\gamma-2j}(V)-c^{-1}_{\gamma-2j}P_{2\gamma-4j}B^{2\gamma}_{2j}(V)\right).
		\end{align*}
		Similarly, for $ \lfloor \gamma/2\rfloor+1\leq j\leq \lfloor \gamma\rfloor$, using \eqref{Pi asy 2} and \eqref{Sec 5 Uj}, we get
		\begin{align*}
			&\rho^{1-n}\partial_{\rho}(\rho^{\frac{n}{2}-\gamma}U_j)\Pi^{\gamma}_{j}(\rho^{\frac{n}{2}-\gamma}(V-\tilde{V}))\big|_{\rho=0}\\
	=&\left(\frac{n}{2}+\gamma-2j\right)F_j\big|_{\rho=0}\cdot\rho^{-\frac{n}{2}+\gamma-2j}\Pi^{\gamma}_{j}(\rho^{\frac{n}{2}-\gamma}(V-\tilde{V}))\big|_{\rho=0}\\
			&+\left(\frac{n}{2}-\gamma+2j\right)G_j\big|_{\rho=0}\cdot\rho^{-\frac{n}{2}-\gamma+2j}\Pi^{\gamma}_{j}(\rho^{\frac{n}{2}-\gamma}(V-\tilde{V}))\big|_{\rho=0}\\
			=&\left(\frac{n}{2}+\gamma-2j\right)b_{2j}B^{2\gamma}_{2\gamma-2j}(U)\left(B^{2\gamma}_{2j}(V)-c^{-1}_{2j-\gamma}P_{4j-2\gamma}B^{2\gamma}_{2\gamma-2j}(V)\right).
		\end{align*}
		For the other term in $A_2$,  for $0\leq j\leq \lfloor \gamma/2\rfloor$,  we have
\begin{align}\nonumber
&\rho^{1-\frac{n}{2}-\gamma}U_j\partial_{\rho}\left[\Pi^{\gamma}_{j}(\rho^{\frac{n}{2}-\gamma}(V-\tilde{V}))\right]\Big|_{\rho=0}\\
\label{4.35}
=&		\rho^{-2j}U_j \Big|_{\rho=0}\cdot \rho^{1-\frac{n}{2}-\gamma+2j}\partial_{\rho}\left[\Pi^{\gamma}_{j}(\rho^{\frac{n}{2}-\gamma}(V-\tilde{V}))\right]\Big|_{\rho=0}\nonumber\\
=&B^{2\gamma}_{2j}(U)\rho^{1-\frac{n}{2}-\gamma+2j}\partial_{\rho}\left[\Pi^{\gamma}_{j}(\rho^{\frac{n}{2}-\gamma}(V-\tilde{V}))\right]\Big|_{\rho=0}.	
\end{align}
By (\ref{Pi asy 1}) and (\ref{Sec 5 formula B 2-2}), we have
		\begin{align}\label{Sec 5 Pi formula a}
&\rho^{1-\frac{n}{2}-\gamma+2j} \partial_{\rho}\left[\Pi^{\gamma}_{j}(\rho^{\frac{n}{2}-\gamma}(V-\tilde{V}))\right]\Big|_{\rho=0}\nonumber\\
			=&\rho^{1-\frac{n}{2}-\gamma+2j}\partial_{\rho}\left[\rho^{\frac{n}{2}+\gamma-2j}\rho^{-\frac{n}{2}-\gamma+2j}
\Pi^{\gamma}_{j}(\rho^{\frac{n}{2}-\gamma}(V-\tilde{V}))\right]\Big|_{\rho=0}\nonumber\\
			=&\left(\frac{n}{2}+\gamma-2j\right)\left[\rho^{-\frac{n}{2}-\gamma+2j}\Pi^{\gamma}_{j}(\rho^{\frac{n}{2}-\gamma}(V-\tilde{V}))\right]\Big|_{\rho=0}+
\nonumber\\
		&\rho\partial_{\rho}\left[\rho^{-\frac{n}{2}-\gamma+2j}\Pi^{\gamma}_{j}(\rho^{\frac{n}{2}-\gamma}(V-\tilde{V}))\right]\Big|_{\rho=0}\nonumber\\	=&\left(\frac{n}{2}+\gamma-2j\right)\left[\rho^{-\frac{n}{2}-\gamma+2j}\Pi^{\gamma}_{j}(\rho^{\frac{n}{2}-\gamma}(V-\tilde{V}))\right]\Big|_{\rho=0}\nonumber\\
=&\left(\frac{n}{2}-\gamma+2j\right)b_{2\gamma-2j}\left(B^{2\gamma}_{2\gamma-2j}(V)-c^{-1}_{\gamma-2j}P_{2\gamma-4j}B^{2\gamma}_{2j}(V)\right).
		\end{align}
Substituting (\ref{Sec 5 Pi formula a}) into (\ref{4.35}), we obtain
		\begin{align}\label{Sec 5 claim 2}
			&\rho^{1-\frac{n}{2}-\gamma}U_j\partial_{\rho}\left[\Pi^{\gamma}_{j}(\rho^{\frac{n}{2}-\gamma}(V-\tilde{V}))\right]\Big|_{\rho=0}\nonumber\\
=&\left(\frac{n}{2}+\gamma-2j\right)b_{2\gamma-2j}B^{2\gamma}_{2j}(U)\left(B^{2\gamma}_{2\gamma-2j}(V)-c^{-1}_{\gamma-2j}P_{2\gamma-4j}B^{2\gamma}_{2j}(V)\right).
		\end{align}
	For $\lfloor \gamma/2 \rfloor + 1 \leq j \leq \lfloor \gamma \rfloor$, a similar argument yields the following identity
		\begin{align*}
		&\rho^{1-\frac{n}{2}-\gamma}U_j\partial_{\rho}\left[\Pi^{\gamma}_{j}(\rho^{\frac{n}{2}-\gamma}(V-\tilde{V}))\right]\Big|_{\rho=0}\\
		=&\left(\frac{n}{2}-\gamma+2j\right)b_{2j}B^{2\gamma}_{2\gamma-2j}(U)\left(B^{2\gamma}_{2j}(V)-c^{-1}_{2j-\gamma}P_{4j-2\gamma}B^{2\gamma}_{2\gamma-2j}(V)\right).
		\end{align*}
		Therefore, we conclude that
		\begin{align*}
			A_2=&2\sum_{j=0}^{\lfloor \gamma/2\rfloor}(\gamma-2j)b_{2\gamma-2j}\int_{\S^n}B^{2\gamma}_{2j}(U)\cdot\left(B^{2\gamma}_{2\gamma-2j}(V)-c^{-1}_{\gamma-2j}P_{2\gamma-4j}B^{2\gamma}_{2j}(V)\right)\ud V_{\S^n}\\
			+&2\sum_{j=\lfloor \gamma/2\rfloor+1}^{\lfloor \gamma\rfloor}(2j-\gamma)b_{2j}\int_{\S^n}B^{2\gamma}_{2\gamma-2j}(U)\left(B^{2\gamma}_{2j}(V)-c^{-1}_{2j-\gamma}P_{4j-2\gamma}B^{2\gamma}_{2\gamma-2j}(V)\right)\ud V_{\S^n}.
		\end{align*}
		The final result follows by defining
		\begin{align*}
		\sigma_{j,\gamma}=&\begin{cases}
		\displaystyle 2(\gamma-2j)b_{2\gamma-2j}, \quad 0\leq j\leq \lfloor \gamma/2\rfloor\\
		\displaystyle 2(2j-\gamma)b_{2j},\quad  \lfloor \gamma/2\rfloor+1\leq j\leq \lfloor \gamma\rfloor
		\end{cases}\\
		=&\begin{cases}
		\displaystyle 2^{2\lfloor \gamma\rfloor+1}j!(\lfloor \gamma\rfloor-j)!\frac{\Gamma(\gamma+1-j)}{\Gamma(\gamma-2j)}\frac{\Gamma(j+1-[\gamma])}{\Gamma(2j+1-\gamma)}, \quad 0\leq j\leq \lfloor \gamma/2\rfloor\\
		\displaystyle -2^{2\lfloor \gamma\rfloor+1}j!(\lfloor \gamma\rfloor-j)!\frac{\Gamma(\gamma+1-j)}{\Gamma(\gamma-2j)}\frac{\Gamma(j+1-[\gamma])}{\Gamma(2j+1-\gamma)},\quad \lfloor \gamma/2\rfloor+1\leq j\leq \lfloor \gamma\rfloor
		\end{cases}
		\end{align*}
		and
		\begin{align*}
	\upzeta_{j,\gamma}=&\begin{cases}
	\displaystyle -2c^{-1}_{\gamma-2j}(\gamma-2j)b_{2\gamma-2j}, \quad 0\leq j\leq \lfloor \gamma/2\rfloor\\
	\displaystyle -2c^{-1}_{2j-\gamma}(2j-\gamma)b_{2j},\quad  \lfloor \gamma/2\rfloor+1\leq j\leq \lfloor \gamma\rfloor
	\end{cases}\\
	=&\begin{cases}
	\displaystyle 2^{4j-2[\gamma]+1}j!(\lfloor \gamma\rfloor-j)!\frac{\Gamma(\gamma+1-j)}{\Gamma(\gamma+1-2j)}\frac{\Gamma(j+1-[\gamma])}{\Gamma(\gamma-2j)}, \quad 0\leq j\leq \lfloor \gamma/2\rfloor\\
	\displaystyle 2^{4\gamma-2[\gamma]+1}j!(\lfloor \gamma\rfloor-j)!\frac{\Gamma(\gamma+1-j)}{\Gamma(2j-\gamma)}\frac{\Gamma(j+1-[\gamma])}{\Gamma(2j+1-\gamma)},\quad \lfloor \gamma/2\rfloor+1\leq j\leq \lfloor \gamma\rfloor.
	\end{cases}
		\end{align*}
	\end{pf}
	
		\textbf{Proof of Theorem \ref{Symmetric Thm}:} To simplify the notation, we define
		\begin{align*}
			&A_2^{'}(U,V)\\
			=&\sum_{j=0}^{\lfloor \gamma/2\rfloor}\upzeta_{j,\gamma}\int_{\S^n}B^{2\gamma}_{2j}UP_{2\gamma-4j}B^{2\gamma}_{2j}V\ud V_{\S^n}+\sum_{j=\lfloor \gamma/2\rfloor+1}^{\lfloor \gamma\rfloor}\upzeta_{j,\gamma}\int_{\S^n}B^{2\gamma}_{2\gamma-2j}UP_{4j-2\gamma}B^{2\gamma}_{2\gamma-2j}V\ud V_{\S^n}.
		\end{align*}
	By the proof of Lemma \ref{Vital integral identity Lem}, we obtain:
	\begin{align*}
		\mathcal{Q}_{2\gamma}(U,V)=A_1(U,V)+A_2^{'}(U,V),
	\end{align*}
	where
	\begin{align*}
		A_1(U,V)=\int_{\B^{n+1}}\rho^{\frac{n}{2}-\gamma}(U-\tilde{U})L^{+}_{2k}(\rho^{\frac{n}{2}-\gamma}(V-\tilde{V}))\ud V_{g_{\B}}.
	\end{align*}
	For the term $ A_1$, it follows from the asymptotic behavior
	\begin{align*}
		\rho^{\frac{n}{2}-\gamma}(U-\tilde{U}), \Delta_{+}^{k}\left(\rho^{\frac{n}{2}-\gamma}(U-\tilde{U})\right)\in  \rho^{\frac{n}{2}-\gamma+2(\lfloor \gamma/2\rfloor+1)}C^{\infty}_{\mathrm{even}}(\overline{\B^{n+1}})+\rho^{\frac{n}{2}+\gamma-2\lfloor \gamma/2\rfloor}C^{\infty}_{\mathrm{even}}(\overline{\B^{n+1}}),
	\end{align*}
where $k\in\mathbb{N}$, and by using Lemma \ref{Green identity}, we conclude that the boundary terms vanish, i.e., $A_1(U, V) = A_1(V, U)$.
By the self-adjointness of $ P_{2\gamma} $ for $ \gamma > 0$, we can then conclude that
	\begin{align*}
	A_2^{'}(U,V)=A_2^{'}(V,U).
	\end{align*}
Thus, we have completed the proof of Theorem \ref{Symmetric Thm}.\\

	\begin{lem}\label{Last lem}
		Let $\gamma\in (0,+\infty)\backslash\mathbb{N}$ and $k\in \mathbb{N}$. For any
$$U\in \rho^{\frac{n}{2}-\gamma+2(\lfloor \gamma/2\rfloor+1)}C^{\infty}_{\mathrm{even}}(\overline{\B^{n+1}})+\rho^{\frac{n}{2}-\gamma+2(\gamma-\lfloor \gamma/2\rfloor)}C^{\infty}_{\mathrm{even}}(\overline{\B^{n+1}}),$$
		we have
		\begin{align}\label{Last lem formula}
		\int_{\B^{n+1}}U(-\tilde{\Delta}_{+})^kU\ud V_{g_{\B}}\geq 0,
		\end{align}
		with  equality if and only if $(-\tilde{\Delta}_{+})^{\lfloor k/2\rfloor}U=0$.
	\end{lem}
	\begin{pf}
	We begin by proving the inequality for $k=1$. In this case, it is equivalent to \begin{align*}
		\int_{\B^{n+1}} U L_{g_{\B}} U  \ud V_{g_{\B}} \geq \frac{1}{4} \int_{\B^{n+1}} U^2  \ud V_{g_{\B}},
	\end{align*}
	where $L_{g_{\B}} = -\Delta_{+} - \frac{n^2-1}{4}$ is the conformal Laplacian on the Poincar\'e ball. Using the conformal invariance of $L_{g_{\B}}$, this reduces to \begin{align}\label{Last lem formula 1}
		\int_{\B^{n+1}} V(x)(-\Delta)V(x) \ud x \geq \int_{\B^{n+1}} \frac{V^2(x)}{(1 - |x|^2)^2}  \ud x,
	\end{align} where
	 \begin{align*}
V(x) = \left(\frac{2}{1 - |x|^2}\right)^{\frac{n-1}{2}} U \in \rho^{\frac{1}{2} - \gamma + 2(\lfloor \gamma/2\rfloor + 1)}C^{\infty}_{\mathrm{even}}(\overline{\B^{n+1}}) + \rho^{\frac{1}{2} + \gamma - 2\lfloor \gamma/2\rfloor}C^{\infty}_{\mathrm{even}}(\overline{\B^{n+1}}).
	 \end{align*}
Notice that, for $\gamma\in (0,\infty)\setminus\mathbb{N}$,
\begin{align*}
2(\lfloor \gamma/2\rfloor + 1)-\gamma>2\cdot \gamma/2-\gamma=0\;\;\;\; \textrm{and}\;\;\;\; \gamma-2\lfloor \gamma/2\rfloor>\gamma-2\cdot \gamma/2=0.
\end{align*}
If we let
 \begin{align*}
\delta_1=2(\lfloor \gamma/2\rfloor + 1)-\gamma\;\;\;\; \textrm{and}\;\;\;\;\delta_2=\gamma-2\lfloor \gamma/2\rfloor,
 \end{align*}
 then $\delta_1>0$, $\delta_2>0$ and
 	 \begin{align*}
V(x)  \in \rho^{\frac{1}{2} +\delta_1}C^{\infty}_{\mathrm{even}}(\overline{\B^{n+1}}) + \rho^{\frac{1}{2} + \delta_2}C^{\infty}_{\mathrm{even}}(\overline{\B^{n+1}}).
	 \end{align*}
Therefore, by integrating by parts, we obtain \begin{align*}
	  \int_{\B^{n+1}} V(x)(-\Delta)V(x)  \ud x = \int_{\B^{n+1}} |\nabla V|^2  \ud x, \end{align*}
	  and inequality \eqref{Last lem formula 1} follows from the Hardy inequality (see e.g. \cite{Brezis&Marcus,Owen}).
	
	Next, for general $k$, observe that for $1 \leq l \leq k$, we have
\begin{align*}
\tilde{\Delta}_{+}^{k-l} U \in \rho^{\frac{n}{2}-\gamma+2(\lfloor \gamma/2\rfloor+1)}C^{\infty}_{\mathrm{even}}(\overline{\B^{n+1}})+ \rho^{\frac{n}{2}-\gamma+2(\gamma-\lfloor \gamma/2\rfloor)}C^{\infty}_{\mathrm{even}}(\overline{\B^{n+1}}),
\end{align*}
which implies
	\begin{align*}
	\rho^{1-n} \tilde{\Delta}_{+}^{l-1} U \times \partial_{\rho} \tilde{\Delta}_{+}^{k-l} U \in \rho^{4(\lfloor \gamma/2\rfloor + 1) - 2\gamma} C^{\infty}_{\mathrm{even}}(\overline{\B^{n+1}}) +
\rho^{2} C^{\infty}_{\mathrm{even}}(\overline{\B^{n+1}})+ \rho^{2\gamma - 4\lfloor \gamma/2\rfloor} C^{\infty}_{\mathrm{even}}(\overline{\B^{n+1}}).
	\end{align*}
	So the boundary terms in Lemma \ref{Green identity} vanish.	That is, if $k = 2m$ with $m \in \mathbb{N}$, then we have
	 \begin{align*} \int_{\B^{n+1}} U(-\tilde{\Delta}_{+})^k U  \ud V{g_{\B}} = \int_{\B^{n+1}} (\tilde{\Delta}_{+}^m U)^2  \ud V{g_{\B}} \geq0;
	 \end{align*}
	  if $k = 2m + 1$, then
	  \begin{align*}
	  	 \int_{\B^{n+1}} U(-\tilde{\Delta}_{+})^k U  \ud V{g_{\B}} = \int_{\B^{n+1}} \left[(-\tilde{\Delta}_{+})^m U\right] (-\tilde{\Delta}_{+}) \left[(-\tilde{\Delta}_{+})^m U\right]  \ud V{g_{\B}} \geq 0,
	  \end{align*}
	  where the final inequality follows from the case $k=1$.
This completes the proof of Lemma \ref{Last lem}.
	\end{pf}

	\textbf{Proof of Theorem \ref{Functional inequ Thm}:} By the definition of   $A_1(U,U)$ and Lemma \ref{Last lem}, we get
	\begin{align*}
		A_1(U,U)=&\int_{\B^{n+1}}	\rho^{\frac{n}{2}-\gamma}(U-\tilde{U})L_{2k}^{+}\left(\rho^{\frac{n}{2}-\gamma}\left(U-\tilde{U}\right)\right)\ud V_{g_{\B}}\\
		=&\int_{\B^{n+1}}	\rho^{\frac{n}{2}-\gamma}(U-\tilde{U})	\prod_{l=0}^{\lfloor \gamma\rfloor}\left(-\tilde{\Delta}_{+}+(\gamma-2l)^2\right)\left(\rho^{\frac{n}{2}-\gamma}\left(U-\tilde{U}\right)\right)\ud V_{g_{\B}}\\
		\geq&\prod_{l=0}^{\lfloor \gamma\rfloor}(\gamma-2l)^2\int_{\B^{n+1}}		\left(\rho^{\frac{n}{2}-\gamma}\left(U-\tilde{U}\right)\right)^2\ud V_{g_{\B}}\geq 0.
	\end{align*}
Therefore,
	\begin{align*}
		\mathcal{Q}_{2\gamma}(U,V)=A_1(U,V)+A_2^{'}(U,V)\geq A_2^{'}(U,V).
	\end{align*}
	The desired result follows by the definition of $A_2^{'}(U,V)$. \\
	
	\textbf{Proof of Corollary \ref{Sobolev Trace Thm}:}  Noticing that, for $\gamma\in (0,\infty)\setminus\mathbb{N}$,
	\begin{align*}
	\gamma\pm\frac{n}{2}\in \mathbb{Z}  \qquad\Longleftrightarrow\qquad[\gamma]=\frac{1}{2}, \quad n\in 2\mathbb{N}-1,
	\end{align*}
where $\mathbb{Z}=\{0,\pm1,\pm2,\cdots\}$ is the set of integers.
	Hence our proof is divide into two parts, we only give the proof of one case.
	
 If $\gamma-\frac{n}{2}\in 2\mathbb{N}\cup\{0\}$, using Theorem \ref{Functional inequ Thm}, there holds
		\begin{align*}
		\mathcal{Q}_{2\gamma}(U,U)\geq &
		\sum_{j=0}^{\lfloor \gamma/2\rfloor}\upzeta_{j,\gamma}\int_{\S^n}B^{2\gamma}_{2j}UP_{2\gamma-4j}B^{2\gamma}_{2j}U\ud V_{\S^n}\\
		+&\sum_{j=\lfloor \gamma/2\rfloor+1}^{\lfloor \gamma\rfloor}\upzeta_{j,\gamma}\int_{\S^n}B^{2\gamma}_{2\gamma-2j}UP_{4j-2\gamma}B^{2\gamma}_{2\gamma-2j}U\ud V_{\S^n}.
		\end{align*}
		By the definition of $j_1$, we have $\gamma-2j_1=\frac{n}{2}$ ,
		\begin{align*}
			\frac{n}{2}<\gamma-2j\in \frac{n}{2}+\mathbb{N}, \qquad\mathrm{for}\qquad 0\leq j<j_1
		\end{align*}
		and
		\begin{align*}
			\gamma-2j<\frac{n}{2} \qquad\mathrm{for}\qquad j_1< j\leq \lfloor \gamma/2\rfloor .
		\end{align*}
		Therefore, using the Theorem A and reverse Sobolev inequalities,   we can get
		\begin{align*}
			&\sum_{j=0}^{\lfloor \gamma/2\rfloor}\upzeta_{j,\gamma}\int_{\S^n}B^{2\gamma}_{2j}UP_{2\gamma-4j}B^{2\gamma}_{2j}U\ud V_{\S^n}\\
			\geq& 2n!|\S^n|\upzeta_{j_1,\gamma}\log \left(\fint_{\S^n} e^{B^{2\gamma}_{2j_1}(U)-\overline{B^{2\gamma}_{2j_1}(U)}}  \ud V_{{\S^n}}\right)\\
			+&\sum_{j=j_1+1}^{\lfloor \gamma/2\rfloor}\frac{\Gamma\left(\frac{n}{2}+\gamma-2j\right)}{\Gamma\left(\frac{n}{2}-\gamma+2j\right)}|\S^n|^{\frac{2\gamma-4j}{n}}
			\upzeta_{j,\gamma}\|B^{2\gamma}_{2j}U\|^2_{L^{\frac{2n}{n-2\gamma+4j}}(\S^n)}.
		\end{align*}
		Since $\gamma-\frac{n}{2}+\gamma+\frac{n}{2}=2\gamma\in 2\mathbb{N}-1$,  we can see $\gamma+\frac{n}{2}\not\in 2\mathbb{N}$, which means  $2j-\gamma\not= \frac{n}{2}$ for $j\geq 1$. Also, there holds $$2(j_2-1)-\gamma<2\left(\frac{\gamma}{2}+\frac{n}{4}\right)-\gamma=\frac{n}{2}<2j_2-\gamma. $$
		We can deduce that
		\begin{align*}
			\frac{n}{2}<2j_2-\gamma=\frac{n}{2}+1<\frac{n}{2}+2.
		\end{align*}
		So we can obtain
			\begin{align*}
		2j-\gamma<\frac{n}{2} \qquad\mathrm{for}\qquad \lfloor \gamma/2\rfloor+1\leq j<j_2
		\end{align*}
		and
		\begin{align*}
	\frac{n}{2}<2j-\gamma\in\frac{n}{2}+\mathbb{N}, \qquad\mathrm{for}\qquad j_2\leq j\leq \lfloor \gamma/2\rfloor .
		\end{align*}
		Then, it is not hard to see
		\begin{align*}
			&\sum_{j=\lfloor \gamma/2\rfloor+1}^{\lfloor \gamma\rfloor}\upzeta_{j,\gamma}\int_{\S^n}B^{2\gamma}_{2\gamma-2j}UP_{4j-2\gamma}B^{2\gamma}_{2\gamma-2j}U\ud V_{\S^n}\\
			\geq& \sum_{j=\lfloor \gamma/2\rfloor+1}^{j_2-1}\frac{\Gamma\left(\frac{n}{2}-\gamma+2j\right)}{\Gamma\left(\frac{n}{2}+\gamma-2j\right)}|\S^n|^{\frac{4j-2\gamma}{n}}
			\upzeta_{j,\gamma}\|B^{2\gamma}_{2\gamma-2j}U\|^2_{L^{\frac{2n}{n+2\gamma-4j}}(\S^n)}.
		\end{align*}
		So  we get the first inequality in the Corollary \ref{Sobolev Trace Thm}. For the proof of the remaining inequalities, we will only highlight the differences relative to the first one. Specifically, this involves comparing
		 $j_1$ and $0$,  $j_2-1$ and  $\lfloor \gamma\rfloor$, and so on.

		\bibliographystyle{unsrt}

	\end{document}